\documentclass[reqno,12pt]{amsart}
\usepackage{mathrsfs}
\usepackage{amssymb,amsmath,amsthm,amsfonts}
\usepackage[dvipsnames]{xcolor}
\usepackage[
		colorlinks=true,
		linkcolor=blue,
		citecolor=ForestGreen,
		urlcolor=cyan
		]
{hyperref}
\usepackage[inline]{enumitem}
\usepackage[margin=1in]{geometry}
\let\oldsquare\square
\usepackage{mathabx}
\usepackage{graphicx}
\usepackage{cases}
\usepackage{cite}
\usepackage{tikz}\usepackage{subfigure}
\usetikzlibrary{decorations.pathreplacing,calligraphy}
\newcommand{\R}{\mathbb R}
\newcommand{\rd}{\R^d}

\newcommand{\Z}{\mathbb Z}

\newcommand{\F}{ {\mathcal F} }
\newcommand{\h}{ {\mathcal H} }

\allowdisplaybreaks

\def\TagOnRight

\theoremstyle{plain}
\newtheorem{theorem}{Theorem} [section]
\newtheorem{lemma}[theorem]{Lemma}
\newtheorem{corollary}[theorem]{Corollary}
\newtheorem{proposition}[theorem]{Proposition}

\theoremstyle{remark}
\newtheorem{remark}[theorem]{Remark}

\theoremstyle{definition}
\newtheorem{definition}[theorem]{Definition}

\def\({\left(}
\def\){\right)}
\def\<{\left\langle}
\def\>{\right\rangle}

\numberwithin{equation}{section}

\begin{document}
\title[Mixed fractional Hartree and Hermite-Hartree equation]
{On the  mixed fractional Hartree and Hermite-Hartree equations  in  modulation and Fourier amalgam spaces}  % in
\author[D. G. Bhimani]{Divyang G. Bhimani}
\address{Department of Mathematics, Indian Institute of Science Education and Research, Dr. Homi Bhabha Road, Pune 411008, India}
\email{divyang.bhimani@iiserpune.ac.in}
%\author{L. Chergui}
\author{Hichem Hajaiej}
\address{Department of Mathematics, College of Natural Science, California State University, 5151 State Drive, Los Angeles, 90032 CA, USA}
\email{hhajaie@calstatela.edu}
\author{Saikatul Haque}
\address{
Department of Mathematics\\
University of California\\Los Angeles\\ 90095 CA\\USA\\
\& Harish-Chandra Research Institute, Allahabad 211019, India}
\email{saikatul@math.ucla.edu}

\thanks{}

\begin{abstract}
We prove local and global well-posedness for mixed fractional Hartree equation and with low regularity Cauchy data in Fourier amalgam $\F W(L^p,\ell^q)$ and modulation  $M^{p,q}$ spaces.   Similar  results also hold for the Hartree equation with  harmonic potential in some modulation spaces.   Our approach also  addresses  Hartree-Fock equations of finitely many (but arbitrary large) particles. A key ingredient of  our method is to  establish trilinear estimates for Hartree non-linearity and the use of Strichartz estimates.   As a consequence,  we could gain $\F W(L^p,\ell^q)$ and $M^{p,q}-$regularity for all $p,q\in [1, \infty].$
In particular, we extend result of  Bhimani-Grillakis-Okoudju \cite{bhimani2020hartree}  in $M^{p,q}$ for all $p,q$ and complement known results in Sobolev spaces.
\end{abstract}
\maketitle
%\tableofcontents
\section{Introduction}
We consider the  Cauchy problem for the  Hartree equation with mixed fractional Laplacian:
\begin{eqnarray}\label{mfh}
\begin{cases}  i\partial_tu - (-\Delta)^{s_1}u - (-\Delta)^{s_2} u=( K \ast |u|^2)u,\\
u(0,x)=u_0(x),
\end{cases} (t, x)\in \mathbb R \times \mathbb R^d.
\end{eqnarray}
Here,  $s_1, s_2 \in \R,  u(t,x),u_0(x)\in \mathbb C$ and $K$
denotes the Hartree kernel
\begin{eqnarray}\label{hk}
K(x)= \tfrac{\lambda}{|x|^{\gamma}}\quad (\lambda \in \mathbb R, 0< \gamma<d, x\in \mathbb R^{d}).
\end{eqnarray}
The classical (i.e. $s_1=s_2=1$ case) and fractional (i.e $s_1=s_2$ case) Hartree equations appear in the several physical  phenomena,  e.g.  optical media,  boson stars,  Brownian motion.     See \cite[Section 1.1]{D3, bhimani2023hartree} and the references therein.
The mixed fractional Laplacian (i.e. $s_1\neq s_2$ case) arises in the
case when a particle can follow
two stochastic processes with a different random walk and a L\'evy flight according to a certain probability. The corresponding limit diffusion is described by a sum of two fractional Laplacians with different orders, see \cite{biagi2022mixed}.
It also models  heart anomalies caused by   arteries issues
by the  superposition of two to five mixed fractional Laplacians, to consider different
anomaly in the five arteries,  see \cite{magin2008modeling}.

In recent  years Cauchy problem for  nonlinear dispersive equations with low regularity  initial data space   have been studied by many authors, see \cite{bhimani2016functions,JH, forlano2020deterministic, bhimani2023hartree, D3, wang2007global, wang2011harmonic, benyi2009local, HerrNA2014}.  In this paper, we establish a local and global well-posedness for \eqref{mfh} with Cauchy data in Fourier  amalgam and modulation spaces.  In order to state our main results, we briefly recall these spaces.
The first appearance of amalgam spaces dates back to the work of  Wiener \cite{NW26,NW} in his study of generalized harmonic analysis, where the  amalgam space $W(L^p, \ell^q)=W(L^p, \ell^q) (\R^d)$ is defined by the norm
\[\|f\|_{W(L^p, \ell^q)}= \left( \sum_{n \in \mathbb Z^d} \left( \int_{n+(0,1]^d} |f(x)|^p dx \right)^{q/p} \right)^{1/q}.\]
In the 1980s, Feichtinger \cite{Fei} introduced  a generalization of amalgam spaces.  This enables a vastly wider range of Banach spaces of  functions  or distributions defined on  locally compact  group to be used as a local or global component, resulting in a deep and powerful theory.  Specifically,  he used the notation $W(B,C)$ to define a space of functions  or distributions which are ``locally in  Banach space $B$" and ``globally in Banach space $C$", and called them \textbf{Wiener amalgam  type} spaces. In order to define these spaces precisely we briefly introduce notations.  For any given function $f$ which is locally in $B$  (i.e,  $gf\in B, \forall g \in C_0^{\infty}(\R^d)),$ we set $f_{B}(x)= \|f g(\cdot -x)\|_{B}$, for some nonzero $g\in C_0^\infty(\rd)$. The space $W(B,C)$ is defined as the space of all functions $f$ locally in $B$ such that  $f_{B}\in C$. The space $W(B, C)$ endowed with the norm  $\|f\|_{W(B, C)}=\|f_{B}\|_{C}.$  Moreover, different choices of nonzero $g\in C^{\infty}_0(\R^d)$  generate the same space and yield equivalent norms, see \cite[Theorem 1]{Fei} and \cite[Proposition 11.3.2]{CH2002}. For an expository introduction to Wiener amalgam spaces on $\mathbb R$  with extensive references to the original literature, we refer to \cite{CH2002, CH2007}.

In this paper we consider the  Fourier image of  a particular Wiener amalgam spaces $W(L^p, \ell^q_s)$,  which is  known as the  \textbf{Fourier amalgam spaces}  $\mathcal{F}W(L^p, \ell^q_s)=\mathcal{F}W(L^p, \ell^q_s)(\R^d).$ 
More specifically, for $1\leq p, q \leq \infty,  s\in \R,$ we  define 
$$\mathcal{F}W(L^p, \ell^q_s)=\bigl\{ f\in \mathcal{S}'(\R^d): \|f\|_{\mathcal{F}W(L^p, \ell^q_s)}= \bigl\|\|\chi_{n+(0,1]^d} (\xi)\mathcal{F}f(\xi)\|_{L_{\xi}^p(\R^d)} \langle n \rangle ^{s}  \bigr\|_{\ell_n^q(\mathbb Z^d)}< \infty \bigr\},
$$
where, $\mathcal{F}$ denotes the Fourier transform,  $\mathcal{S}'(\R^d)$ is the space of tempered distributions.

 Now we turn our attention to modulation spaces, which were born during the early eighties in pioneering work of H. Feichtinger \cite{Fei}. It  is now present in both pure and applied mathematics and appeared in many applications, see e.g. \cite{HGF2006, HF2015, kassob, FeiPAMS,grochenig2013foundations}. In particular, we note that it has played a central role in the long standing quest  to  understand the dispersive PDEs (e.g  NLS and mKdV) near scaling criticality in the last two decades, see \cite{ruzhansky2012modulation, wang2011harmonic, wang2007global}.
In contrast with the Besov spaces, which are defined by a dyadic decomposition of the frequency space, modulation spaces  arise from a uniform partition of the frequency space. In order to  make this definition precise, we  introduce some notations.  Let us start with uniform covering of $\mathbb R^n$ by unit cubes, specifically,  $\mathbb R^d = \bigcup_{k\in \mathbb Z^{d}} Q_{k},$ where $Q_{k} \subset \mathbb R^d$ be  the unit cube with center at $k.$  Let   $\rho \in \mathcal{S}(\R^d)$ (Schwartz space),  $\rho: \R^d \to [0,1]$  be  a smooth function satisfying   $\rho(\xi)= 1 \  \text{if} \ \ |\xi|_{\infty}\leq \frac{1}{2} $ \footnote{Define $|\xi|_{\infty}=\max\{ | \xi_i | : \xi= (\xi_1,..., \xi_n)\}.$} and $\rho(\xi)=
0 \  \text{if} \ \ |\xi|_{\infty}\geq  1.$ Let  $\rho_n$ be a translation of $\rho,$ that is,
$ \rho_n(\xi)= \rho(\xi -n), n \in \Z^d$
and denote
$\sigma_{n}(\xi)=
  \frac{\rho_{n}(\xi)}{\sum_{\ell\in\Z^{d}}\rho_{\ell}(\xi)},  n
  \in \Z^d.$ This family of smooth functions gives  a   bounded admissible partition of unity  (BAPU), i.e. $\{\rho_n\}_{n \in \mathbb Z^d}$ satisfies the following conditions:
\begin{align*}
    \begin{cases}
    |\sigma_{k}(\xi)|\geq c, \forall z \in Q_{k}, \text{ for some } c>0\\
    \text{supp}   \sigma_{k} \subset \{\xi: |\xi-k|_{\infty}\leq 1 \}\\
\sum_{k\in \mathbb Z^{d}} \sigma_{k}(\xi)\equiv 1, \quad  \forall \xi \in \mathbb R^d\\
  |D^{\alpha}\sigma_{k}(\xi)|\leq C_{|\alpha|}, \forall \xi \in \mathbb R^d, \quad  \alpha \in (\mathbb N \cup \{0\})^{d}.
    \end{cases}
\end{align*}
Then the frequency-uniform decomposition operators can be defined by
\[\oldsquare_n = \mathcal{F}^{-1} \sigma_n \mathcal{F}. \]
The \textbf{modulation  space}  $M^{p,q}_s(\R^d)$ is defined by the norms:
\begin{equation*}
\|f\|_{M^{p,q}_s}=   \bigl\| \left\lVert
  \oldsquare_nf\right\rVert_{L_x^p(\R^d)} \langle n \rangle ^{s} \bigr\|_{\ell^q_n(\mathbb Z^d)}.
\end{equation*} For $s=0,$ we write $M^{p,q}_0(\rd)=M^{p,q}(\rd)$ and $\F W(L^p,\ell_0^q)=\F W(L^p,\ell^q)$,  see also Remark \ref{stft}. %and \ref{chr} \textcolor{red}{(erase it?)}.  
The \textbf{Fourier-Lebesgue spaces} $\mathcal{F}L_s^p(\R^d)$ is defined by
$$\mathcal{F}L_s^p(\R^d)= \bigl\{f\in \mathcal{S}'(\R^d): \|f\|_{\mathcal{F}L_s^{p}}:= \|\langle\xi\rangle^s\F{f}\|_{L^p}< \infty \bigr\}.$$  It turns out that
\begin{eqnarray*}
\F W(L^p,\ell_s^q)=
\begin{cases}M_s^{2,q}\ (\text{modulation spaces}) \quad& \text{if}\ p=2\\
\mathcal{F}L_s^q \ ( \text{Fourier-Lebesgue spaces}) \quad &\text{if} \  p=q\\
H^s \ (\text{Sobolev space}) \quad &\text{if} \  p=q=2.
\end{cases}
\end{eqnarray*}
We refer to \cite{wang2011harmonic, ruzhansky2012modulation, kassob, forlano2020deterministic, CH2002, CH2007, HGF2006,feichtinger1983modulation,HF2015,grochenig2013foundations} for comprehensive introduction to these spaces. Let us denote  $Y_{rad}$  the space of radial functions in $Y$.
\begin{theorem}[local well-posedness]\label{MT} Let $K$ be given by \eqref{hk} and
\begin{eqnarray*}
X=\begin{cases} \F W(L^p,\ell^q) & \text{if}  \    1 \leq q \leq \frac{2d}{d+\gamma}\leq p \leq \infty, s_1, s_2 \in \R, 0<\gamma<d\\
\F W(L^p,\ell^q)\cap L^2 &  \text{if} \    1\leq q<\frac{2d}{d+\gamma}, q\leq p\leq\infty, s_1, s_2 \in \R, 0< \gamma<d \\
\F W(L^p,\ell^q)\cap L^2 &  \text{if}   \ p,q \in [1, \infty],  d\geq 1, s_1=s_2=1,  0<\gamma<\min(2,\frac{d}{2})\\
\F W(L^p,\ell^q)\cap L_{rad}^2 &  \text{if}   \  p,q \in [1, \infty],  d\geq 2,  \frac{d}{2d-1}< s_1\leq s_2 \leq 1, 0<\gamma<\min(2s_2,\frac{d}{2})
\end{cases}
\end{eqnarray*}
or
\begin{eqnarray*}
X=\begin{cases}
M^{p,q}\cap L^2 &  \text{if}   \ p,q \in [1, \infty],  d\geq 1, s_1=s_2=1,  \gamma<\min(2,\frac{d}{2})\\
M^{p,q}\cap L_{rad}^2 &  \text{if}   \ p,q \in [1, \infty],  d\geq 2,  \frac{d}{2d-1}< s_1\leq s_2 \leq 1, \gamma<\min(2s_2,\frac{d}{2}).
\end{cases}
\end{eqnarray*}
Then \eqref{mfh} is locally well-posed in $X$:
for given $m>0$, there exist $T=T(m)>0$ and $X_T\subset  C([0,T),X)$
such that for each $u_0\in B_m=\{f\in X:\|f\|_X\leq m\}$, \eqref{mfh} has a unique solution in $X_T$.
Moreover the solution map $u_0\mapsto u$ is (Lipschitz) continuous from $B_m$ to $C([0,T),X)$.
\end{theorem}
Theorem \ref{MT} is new  even for $s_1=s_2=1$ (classical Hartree equation) in Fourier amalgam spaces $\F W(L^p,\ell^q)$ for $p\neq q.$   Note that there is $f\in  \F W(L^p,\ell^q)$ (with $1 \leq q < \frac{2d}{d+\gamma}  < p < 2$) such that $f\notin L^2$ (see Lemma \ref{pl} \eqref{exp}).  Thus,  in particular,  we could establish local well-posedness for some initial data with infinite $L^2-$norm.   We would like to mention that Vargas and Vega in \cite{vargas2001JMPA} have   studied  the well-posedness for 1D cubic NLS for initial data with infinite $L^2-$norm.  The idea of their  proof is to  decompose the initial data  into two parts that  satisfy certain suitable estimates. This strategy was originally introduced by J. Bourgain in \cite{bourgain1999JAMS}, see also \cite{leonid2017JDE,  hyakuna2012PAMS,manna2022JDE}.
On the other hand, our approach is rather simple, but we could still cover rough initial data (which was not covered earlier for the Hartree equation). The key ingredients in our proof is the  fact that any unimodular   Fourier multiplier operator is unitary on  Fourier amalgam spaces (Lemma \ref{pl} \eqref{n1}),  new  trilinear  estimates (Propositions  \ref{t2} and \ref{t4}) and Strichartz estimates. 

We note that  Bhimani et al.  in   \cite[Theorem 1.1]{D3} established local well-posedness for \eqref{hH} in $M^{p,q} \subset L^2$ for $1\leq p \leq 2,  1\leq q \leq \frac{2d}{d+\gamma}.$ Their approach was based on  trilinear estimates  and boundedness of  Fourier multiplier in $M^{p,q}.$ The novelty in this paper is the  use of Strichartz estimates (Theorem \ref{stri}) in addition to %Bhimani et al.  
ideas in \cite{D3},    to obtain the full range of $p,q\in [1, \infty].$

\begin{remark}  Theorem \ref{MT} deserves several further comments.
\begin{enumerate}

\item Herr-Lenzmann \cite[Theorem 1.1]{HerrNA2014} established  local well-posedness for boson star equation  $i \partial_tu - (-\Delta + m^2)^{1/2}+ (|x|^{-1} \ast |u|^2)u =0,$ where $m\geq 0,$   in  $H^s{(\R^3)}$ for $s>1/4$ and in $H^s_{rad}$ for $s>0$.    In view of Lemma \ref{pl}\eqref{n1},  we note that the first cases of Theorem \ref{MT} hold for the boson star equation.   Since $H^{s}\subset  \mathcal{F}W(L^2, \ell^1) \subset L^2$  for any $s>0,$ Theorem \ref{MT} complements \cite[Theorem 1.1]{HerrNA2014} as we do not need any radial assumption  for initial data and it works for all dimensions.

\item    The single fractional  Hartree equation,  i.e.  \eqref{mfh}  with $s_1=s_2 \in (1/2,1] $ is locally well-posed in $H^{s}$ for $s\geq \frac{\gamma}{2}$ and in $\mathcal{F}L^q$  ($q\in [2d/(d+\gamma), 2]$) for $s_1=s_2=1, 0<\gamma<\min (2,d)$ or $\gamma=2, d\geq 3;$ see \cite{hichem2013FE, changxing2008cauchy, HyakunaMulti, cho2017AA, hajaiej2014multilinear,  hajaiej2011necessary}.  Cho-Ozawa \cite{cho2006SIAM}  established  several local well-posedness  for boson star equation.
Since $\mathcal{F}W(L^q, \ell^q) =\mathcal{F}L^q$ and in view of  embedding between $\F W(L^p,\ell^q)$ and $H^s,$  Theorem \ref{MT}
complements these results.
\item
In \cite[Theorem 1.5]{bhimani2023hartree},  Bhimani-Haque proved local well-posedness for Hartree equation in  $\mathcal{F}L^q\cap L^2$ for $1\leq q \leq \infty$ by employing Strichartz estimates. Earlier only the case $q\geq2$ was known. Theorem \ref{MT} recovers this result and extends to mixed fractional Laplacian.
\end{enumerate}
\end{remark}
\begin{remark} We discuss several imposed  hypotheses of Theorem \ref{MT}.
\begin{enumerate}

\item In view of  Lemma \ref{pl} \eqref{n1},   the trilinear estimate in  Proposition \ref{t2} gives local well-posedness for data in $\F W(L^p,\ell^q)$ ($1 \leq q \leq \frac{2d}{d+\gamma}\leq p \leq \infty$).  We shall notice that these restrictions on $p,q$ comes due to Hardy-Littlewood Sobolev  inequality and inclusion  relation of these spaces in Lemma \ref{pl} \eqref{in}.

\item  Note that we can take $p< \frac{2d}{d+\gamma}$ for data in $\F W(L^p,\ell^q)\cap L^2$ ($1\leq q<\frac{2d}{d+\gamma}, q\leq p\leq\infty$).   Here the restriction on $p, q$ comes due to Proposition \ref{t4} \eqref{s1}.   In fact,  we shall need this in order to insure   that $k_2$ (the second part of the Fourier transform the Hartree kernel,  see \eqref{k1k2}) is in $L^{q/(2(q-1))}.$

\item Taking dispersion   $s_1=s_2=1$, i.e. classical Laplacian  $(-\Delta)$ in \eqref{mfh},  and $  0<\gamma<\min(2,\frac{d}{2}),$  we could  employ Strichartz estimate (Theorem \ref{stri} \eqref{stri2}),  and we do not require any restriction on $p,q,$ i.e.   $1\leq p, q \leq \infty$ for data in $\F W(L^p,\ell^q)\cap L_{rad}^2$ or in $M^{p,q}\cap L^2$.  We shall notice that  in order to use Strichartz estimate we impose $\gamma<2$ and to use Proposition \ref{t4}, we impose $\gamma<\frac{d}{2}$.

\item  In order to treat  mixed fractional Laplacian $(\frac{d}{2d-1}< s_1\leq s_2 \leq 1)$,   we need to  impose radiality assumption on initial data,  dimension   restriction $d\geq 2$    due to Strichartz estimate Theorem \ref{stri} \eqref{stri1}.
\item Consider the nonlinear Schr\"odinger equation (NLS)  with local nonlinearity:
$$iu_t + \Delta u \pm |u|^2u =0.$$
In    \cite[Theorem 1.1]{benyi2009local} and in \cite[Theorem 4.1]{bhimani2016functions}, authors have established local well-posedness  for NLS in weighted modulation spaces     $M_s^{p,q}$ ($1\leq p, q \leq \infty$) for $s> d (1-1/q).$  See also  \cite[Theorem 1.1]{wangJFA2006} and  \cite[Theorem 1]{klausJFAA2023}. While in \cite{forlano2020deterministic}  1D cubic NLS is studied in the Fourier amalgam spaces.
Compared with   with  local   non-linearity $|u|^{2}u$, the \textit{non-local}  nonlinearity $(|\cdot|^{-\gamma} \ast |u|^2)u$    bring some new   difficulties in order to prove trilinear estimates.   Apparently local well-posedness,  even for classical  Hartree equation,  remain open  in weighted modulation and in weighted  Fourier amalgam spaces. 
\end{enumerate}
\end{remark}
\begin{theorem}[global  well-posedness]\label{mt2}
Let  $K$ be given by \eqref{hk} with  $0 < \gamma <  \min\{2s_2, d/2\}$, $0<s_1\leq s_2\leq1$  and
\begin{eqnarray*}
X=\begin{cases}
\F W(L^p,\ell^q)\cap L^2 &  \text{if} \  \  1\leq q<\frac{2d}{d+\gamma}, q\leq p\leq\infty,  s_1=s_2=1, d\geq 1 \\
\F W(L^p,\ell^q)\cap L_{rad}^2 &  \text{if}   \ p,q \in [1, \infty],  d\geq 2,  \frac{d}{2d-1}< s_1\leq s_2 \leq 1, \gamma<2s_2.
\end{cases}
\end{eqnarray*}
or
\begin{eqnarray*}
X=\begin{cases}
M^{p,q}\cap L^2 &  \text{if}   \ p,q \in [1, \infty],  d\geq 1, s_1=s_2=1,  \gamma<2\\
M^{p,q}\cap L_{rad}^2 &  \text{if}   \ p,q \in [1, \infty],  d\geq 2,  \frac{d}{2d-1}< s_1\leq s_2 \leq 1, \gamma<2s_2.
\end{cases}
\end{eqnarray*}
Assume that $u_0\in X.$ Then there exists a unique global solution $u$ of \eqref{mfh}  such that
$$u \in C(\mathbb R, X)
\cap L^{8s_j/\gamma}_{loc}(\mathbb R, L^{4d/(2d-\gamma)}) .$$
\end{theorem}
Theorem \ref{mt2} is the first global well-posedness result for mixed fractional Hartree equation  \eqref{hH} as far as we are aware.  In \cite{hichem2013FE},  the second author jointly  with Y. Cho,  G. Hwang  and T. Ozawa  have studied  well-posed for fractional Hartree equation.  Theorem \ref{mt2} complements this work.   Carles-Mouzaoui  \cite{carles2014cauchy} proved global well-posedness for classical Hartree equation in the  Wiener algebra $\mathcal{F}L^1 \cap L^2.$  The idea of their proof is to use  global well-posedness results in $L^2$,  conservation of mass (see Proposition \ref{miD}),  and
the algebraic properties of Wiener algebra,  specifically it is a $\mathcal{F}L^1-$module (see Corollary \ref{m1}),  see also \cite{D3, bhimani2023hartree}.   We employ this  strategy in Theorem \ref{mt2} in order to get the Fourier amalgam and modulation space regularity.
Of course,  there is an extensive literature (e.g.  \cite{HerrNA2014, changxing2008cauchy, hichem2013FE, GuoJDE})  on the Hartree equation that assures $H^s-$regularity, Theorem \ref{mt2} assures  $\F W(L^p,\ell^q)$ and $M^{p,q}-$ regularity.
\subsection{Hartree equation with harmonic potential}
The Hermite operator (also known as quantum harmonic oscillator) $H=-\Delta+|x|^2$ plays a vital role in quantum mechanics and analysis (see e.g. \cite{D3, bhimani2018nonlinear} the references therein).
 The  spectral decomposition  of $H$ on $\mathbb R^d$ is given by
\begin{equation}\label{eq spec dec}
H=\sum_{k=0}^{+\infty}(2k+d)P_k,\qquad P_k f=\sum_{|\alpha|=k}\langle f,\Phi_\alpha\rangle\Phi_\alpha,
\end{equation}
where $\langle\cdot,\cdot\rangle$ is the inner product in $L^2$ and $\Phi_\alpha$, $\alpha\in \mathbb{N}^{d}$, are the normalised Hermite functions, forming an orthonormal basis for  $L^2$.
We remark  that $ H^\beta $ is a densely defined unbounded operator.   We thus define the Schr\"odinger  semigroup associated to $H$  by
\begin{equation}\label{shH}
e^{-itH}f=\sum_{k=0}^{+\infty} e^{-it(2k+d)} P_kf.
\end{equation}

We now  consider  Hartree equation with harmonic potential $H=-\Delta+|x|^2$:
\begin{eqnarray}\label{hH}
\begin{cases}  i\partial_tu - (-\Delta+ |x|^2)u =( K \ast |u|^2)u,\\
u(0,x)=u_0(x),
\end{cases} (t, x)\in \mathbb R \times \mathbb R^d.
\end{eqnarray}
In this setting,  we have the following theorem.
 \begin{theorem}\label{gHh} Let $1\leq p <\infty.$ Then
 \begin{enumerate}
 \item \label{gHh1} \eqref{hH} is locally well-posed in $M^{p,p}$ for $0<\gamma <\min\{2,d/2\}$.
 \item\label{gH2} \eqref{hH} is globally well-posed in $M^{p,p} \cap L^2$ for $0<\gamma <\min \{2, d/2\}:$ Given $u_0\in M^{p,p} \cap L^2,$  there exists a unique global solution $u$ of \eqref{hH}  such that
$$u \in C(\mathbb R, M^{p,p})
\cap L^{8/\gamma}_{loc}(\mathbb R, L^{4d/(2d-\gamma)} ).$$
 \end{enumerate}
 \end{theorem}
Bhimani et al.  in \cite[Theorem 1.3 ]{D3} proved global well-posedness for \eqref{hH} in $M^{p,p}$ for $1\leq p \leq \frac{2d}{d+\gamma}$.  Theorem \ref{gHh} extends this result for all $p\in [1, \infty).$ See also  \cite[Theorem 1.1]{bhimani2018nonlinear}. 
\begin{remark} We have the following comments for Theorem \ref{gHh}.
\begin{enumerate}
\item   We consider data in $M^{p,  p}$ mainly because   $e^{-itH}$ (see \eqref{shH}) is  bounded on $M^{p, p}.$ See Proposition \ref{mso} below.   We note that $e^{-itH}$ may not be bounded on $M^{p,q}$  for $p\neq q,$ see  \cite[Remark 4]{bhimani2018nonlinear}.   In view of this,  we cannot even expect to solve  free  \eqref{hH},  i.e.  \eqref{hH} without non-linearity,   in  $M^{p,q}$ for $p\neq q$. 

\item  Recently Bhimani et al.  in \cite{bhimaniAdM2021, bhimani2022heat} have  carried out some interesting study for  heat equation associated to fractional harmonic oscillator $H^{\beta} (\beta>0)$ in modulation spaces.
 It remains interesting  open question  to study \eqref{hH} associated to  $H^\beta$.
\end{enumerate}
\end{remark}

\subsection{Further Remarks} % Hartree not Hartre,  hgfei 
\begin{remark} In \cite[Theorem 1.2]{D3},  the first author jointly with M. Grillakis and K.  Okoudjou have established  global well-posedness for Hartree-Fock equation of finite particles with single fractional Laplacian in some modulation spaces.  Exploiting the ideas of the proof of Theorem \ref{mt2},  this result (i.e.  \cite[Theorem 1.2 ]{D3})  can be generalized to mixed fractional Laplacian.
\end{remark}
\begin{remark}The analogue of Theorem \ref{mt2} also holds true for  reduced Hartree-Fock and Hartree-Fock equations of finitely many particles in Fourier amalgam spaces.  For presenting the clarity of ideas,  instead we have chosen to discuss the single particle equation.
\end{remark}
\begin{remark}
The short-time Fourier transform (STFT) of a  $f\in \mathcal{S}'(\R^d)$ with
 respect to a nonzero window function $ g \in {\mathcal S}(\R^d)$ is defined by
\begin{equation*}\label{stft}
V_{g}f(x,y)= \int_{\R^d} f(t) \overline{g(t-x)} e^{- 2\pi i y\cdot t}dt,  \  (x, y) \in \R^d \times \R^d
\end{equation*}
 whenever the integral exists.   It is known \cite[Proposition 2.1]{wang2007global}, \cite{feichtinger1983modulation} that
\[ \|f\|_{M^{p,q}_s}\asymp  \bigl\| \|V_gf(x,y)\|_{L^p(\R^d)} \langle y \rangle^s \bigr\|_{L^q(\R^d)} \]
 The definition of the modulation space  is independent of the choice of
the particular window function, see e.g. \cite[Proposition 11.3.2(c)]{grochenig2013foundations} and \cite{feichtinger1983modulation}.
 \end{remark}

\begin{remark} Quantative well-posedness is a stronger version of well-posedness, we refer to \cite[Section 3]{bejenaru2006sharp} for definition of this abstract notion.  In view of the uniform boundedness of Schr\"odinger propagator  (Lemma \ref{pl} \eqref{n1}) and trilinear estimates (Proposition \ref{t2} and \ref{t4}\eqref{s1}),  for the first two case of Threorem \ref{MT} i.e. $X=\F W(L^p,\ell^q)$, $q\in[1,\frac{2d}{d+\gamma}],p\in[\frac{2d}{d+\gamma},\infty]$ or $\F W(L^p,\ell^q)\cap L^2$, $q\in[1,\frac{2d}{d+\gamma}),p\in[q,\infty]$, we have \eqref{mfh} is quantitatively well posed (and hence analytically well-posed) in $X$, $C([0,T),X)$. % analyticity of solution in the first two cases for $X$,
In this case for small enough $\|u_0\|_X$,  one can write the solution $u$ as a power series expansion
\[
u=\sum_{k=1}^\infty A_k[u_0]
\]where $A_k[u_0]$ is certain   sum over $3-$ary trees with $k$ nodes, see
\cite[Theorem 3]{bejenaru2006sharp}.
\end{remark}
\begin{remark}
In \cite{bhimani2023sharp}, mixed fractional NLS with inhomogeneous nonlinearity is treated in $L^2-$based Sobolev spaces by establishing Strichartz estimates in Lorentz spaces.
 \end{remark}
\section{Key Estimates}
\subsection{Preliminaries}
   The notation $A \lesssim B $ means $A \leq cB$ for  some constant $c > 0 $ independent of $A,B$. The symbol $\langle k\rangle$ stands for $\sqrt{1+|k|^2}$ for $k\in\rd$.
We recall convolution property for amalgam type spaces:
\begin{theorem}[see Theorem 11.8.3 in \cite{CH2002} and \cite{Fei}]\label{fs} If $L^{p_1}\ast L^{p_2}\subset L^{p_3}$ and $\ell^{q_1}\ast \ell^{q_2}\subset \ell^{q}$ then 
    $W(L^{p_1}, \ell^{q_1}) \ast W(L^{p_2}, \ell^{q_2}) \subset W(L^{p}, \ell^{q})$. Moreover, there is a constant $C>0$ such that for all $f \in W(L^{p_1}, \ell^{q_1})$ and $g \in W(L^{p_2}, \ell^{q_2})$, we have 
    \[\|f\ast g\|_{W(L^{p}, \ell^q)} \leq C \|f\|_{W(L^{p_1}, \ell^{1})} \| g\|_{W(L^{p_2}, \ell^{2})}.\]
\end{theorem}
\begin{corollary}[Pointwise multiplication]\label{m1} Let $\frac{1}{p_1}+ \frac{1}{p_2}=1+ \frac{1}{p}$ and $\frac{1}{q_1}+\frac{1}{q_2}=1+\frac{1}{q}.$ There is a constant $C>0$ such that for all $f \in \mathcal{F}W(L^{p_1}, \ell^{q_1})$ and $g \in \mathcal{F}W(L^{p_2}, \ell^{q_2})$, we have 
\[\|fg\|_{\mathcal{F}W(L^p, \ell^q)}\leq C \|f\|_{\mathcal{F}W(L^{p_1}, \ell^{q_1})}\|g\|_{\mathcal{F}W{(L^{p_2},\ell^{q_2}})}.\]
In particular, 
$\mathcal{F}W(L^p, \ell^q)$ is an $\F L^1$-module i.e. $\|fg\|_{\mathcal{F}W(L^p, \ell^q)}\lesssim\|f\|_{\F L^1}\|g\|_{\mathcal{F}W(L^p, \ell^q)}$.
\end{corollary}
\begin{proof} Recall that  Young's convolution inequality gives $$L^{p_1}\ast L^{p_2}\subset L^{p} \quad \text {and} \quad \ell^{q_1}\ast \ell^{q_2}\subset \ell^{q}.$$ 
Taking this into account,  Theorem \ref{fs} gives
\begin{align*}
   \|fg\|_{\mathcal{F}W(L^p, \ell^q)} & = \|\widehat{f}\ast \widehat{g}\|_{W(L^p, \ell^q)} \\
   & \leq C \|\widehat{f}\|_{W(L^{p_1}, \ell^{q_1})}\|\widehat{g}\|_{W(L^{p_2}, \ell^{q_2})}\\
   & =  C \|f\|_{\mathcal{F}W(L^{p_1}, \ell^{q_1})}\|g\|_{\mathcal{F}W{(L^{p_2},\ell^{q_2}})}.\qedhere
\end{align*}    
\end{proof}   
\begin{lemma}[Basic properties in $\F W(L^p,\ell^q)$,  see e.g.  Lemma 2.1 in \cite{bhimani2023nodea}, \cite{CH2002, CH2007,forlano2020deterministic}] \label{pl} Let $p_j, q_j, p, q \in [1, \infty] $ and $ s_i,  s \in \R$,   where $j=1,2$.
\begin{enumerate}
\item \label{in}(inclusion) $\F W(L^{p_1},\ell^{q_1})\hookrightarrow\F W(L^{p_2},\ell^{q_2})$ for $p_1\geq p_2, q_1 \leq q_2.$
%\item \label{pm} (pointwise multiplication)  Let $\frac{1}{p_1}+ \frac{1}{p_2}=1+ \frac{1}{p}$ and $\frac{1}{q_1}+\frac{1}{q_2}=1+\frac{1}{q}.$ Then we have\[\|fg\|_{\F W(L^p,\ell^q)}\lesssim \|f\|_{\widehat{w}^{p_1,q_1}}\|g\|_{\widehat{w}^{p_2,q_2}}.\]In particular,$\F W(L^p,\ell^q)$ is an $\F L^1$-module i.e. $\|fg\|_{\F W(L^p,\ell^q)}\lesssim\|f\|_{\F L^1}\|g\|_{\F W(L^p,\ell^q)}$.
\item \label{con} (convolution inequality) $\|f\ast g\|_{\F W(L^p,\ell^q)}\leq\|f\|_{\F L^\infty}\| g\|_{\F W(L^p,\ell^q)}$.
\item\label{n1} (uniform boundedness  of linear propogator) Define the Schr\"odinger propogator associated to mixed fractional Laplacian by
\begin{equation}\label{0}
\F U(t)f(\xi)=e^{it(|\xi|^{2s_1}+|\xi|^{2s_2})}\F f.
\end{equation}
Then \[
\|U(t)f\|_{\F W(L^p,\ell_s^q)}=\|f\|_{\F W(L^p,\ell_s^q)}.
\]
In fact,   the same estimate hold for any  $\sigma (\xi)$ (real function) symbol in the Fourier space,  i.e. for  $ \F U(t)f(\xi)=e^{it \sigma (\xi)}\F f.$
\item \label{exp} (examples) (i) There is $f\in \F W(L^p,\ell^q)$  ($1\leq p<2, 1\leq q \leq \infty$) that is not in $L^2$. (ii) $\F W(L^p,\ell^q)$ ($1\leq p\leq 2 \leq q \leq \infty$)  are larger spaces than  $L^2.$ In fact,   we have 
$$L^2\subset M^{2, q}=\F W(L^2,\ell^q) \subset \F W(L^p,\ell^q)\subset\F W(L^1,\ell^\infty)$$ In particular,   $\F W(L^1,\ell^\infty)$ is the largest space in  these family of spaces.
\end{enumerate}
\begin{proof}We only write the proof for \eqref{exp}, as others' proof can be found in the references mentioned above.
Let $f$ be given by $\widehat{f}(\xi)=\chi_{\{|\xi|\leq1\}}\frac{1}{|\xi|^{d/2}}$ then clearly $f\not\in L^2$. But for $|n|\geq2$, $\|\chi_{n+Q}(\xi)\widehat{f}(\xi)\|_{L_\xi^p}=0$ and for $1\leq p<2$, $\|\chi_{n+Q}(\xi)\widehat{f}(\xi)\|_{L_\xi^p}<\infty$. Hence for any $1\leq p<2$, $1\leq q\leq\infty$, we get $f\in\F W(L^p,\ell^q)$.
\end{proof}
\end{lemma}
\begin{proposition}[Basic properties in $M^{p,q}$ \cite{wang2011harmonic, kassob}]\label{des}   Let $s\in \mathbb R$ and $1\leq p, q \leq \infty.$ Define
$$ e^{it\phi (h(D))}f(x)=  \int_{\mathbb R^d}  e^{i \pi t \phi \circ h(\xi)}\, \widehat{f}(\xi) \, e^{2\pi i \xi \cdot x} \, d\xi$$ for $f\in \mathcal{S}(\mathbb R^d)$,
where  $\phi\circ h:\mathbb R^d \to \mathbb R$ is the composition function of $h:\mathbb R^d\to \mathbb R$ and $\phi:\mathbb R \to \mathbb R.$   Let $1\leq p,q\leq\infty$ and $s\in\R.$
\begin{enumerate}
\item (\cite[Theorem 1.1]{deng2013estimate}) \label{ds0}  Assume that   there exist $m_1, m_2>0$ such that  $\phi$ satisfies

$$\left| \phi^{(\mu)} (r)\right| \lesssim\begin{cases}
  r^{m_1-\mu} \quad \textrm{if}\quad  r\geq 1\\
  r^{m_2-\mu}  \quad  \textrm{if}\quad  0<r<1\end{cases}$$
for all  $\mu \in \mathbb N_0$  and $h\in C^{\infty}(\mathbb R^d \setminus \{ 0\})$ is a positive homogeneous function  with degree $\lambda>0$. Then
\[\left\| e^{it\phi (h (D))} f \right\|_{M^{p,q}_s} \lesssim  \|f\|_{M^{p,q}_s}  + |t|^{d \left| \frac{1}{2}- \frac{1}{p} \right|} \|f\|_{M^{p,q}_{s+ \gamma (m_1, \lambda)}}\]
where $\gamma (m_1, \lambda)=d (m_1\lambda -2) |1/2-1/p|$.
\item\label{des2} (consequence of above part \eqref{ds0} and inclusion relation in modulation spaces) Let $0<s_1\leq s_2\leq1$ and $U(t)$ be as in \eqref{0}.
Then
\[
\|U(t)f\|_{M_s^{p,q}}\lesssim(1+|t|^{d \left| \frac{1}{2}- \frac{1}{p} \right|})\|f\|_{M^{p,q}_s}.
\]
\item \label{mfl} $M^{p,q}$ is an $\F L^1$-module i.e. $\|fg\|_{M^{p,q}}\lesssim\|f\|_{\F L^1}\|g\|_{M^{p,q}}$.
\end{enumerate}
\end{proposition}
\subsection{Trilinear Estimates}
We denote  Hartree nonlinearity  by
\begin{eqnarray*}
H_{\gamma}(f, g, h):= \left( K \ast (f \bar{g}) \right) h \quad  (f, g, h \in \mathcal{S}(\mathbb R^d)),
\end{eqnarray*}
where $K$ is given by \eqref{hk}.
The Fourier transform of $K$  is given by
\begin{eqnarray*}
\widehat{K}(\xi)= \tfrac{c}{|\xi|^{d-\gamma}},
\end{eqnarray*}
where $c=\lambda C(d, \gamma)$ is a  constant.
Note that  $\widehat{K}$  does not belong to $L^p-$spaces.  However, we can decompose $\widehat{K}$ into Lebesgue spaces:
\[\widehat{K}=k_1+k_2, \]
where
\begin{equation}\label{k1k2}
\begin{cases} k_{1}:= c\chi_{\{|\xi|\leq 1\}}|\cdot|^{\gamma-d} \in L^r(\mathbb R^{d})\quad \forall \ r\in [1, \frac{d}{d-\gamma})\\
k_{2}:=c \chi_{\{|\xi|>1\}} |\cdot|^{\gamma-d}\in L^r(\mathbb R^{d}) \quad \forall\ r \in (\tfrac{d}{d-\gamma}, \infty].
\end{cases}
\end{equation}
\begin{proposition} \label{t2}  Let $0<\gamma <d,$  $\frac{2d}{d+\gamma} \leq p \leq \infty$ and $ 1\leq q \leq \frac{2d}{d+\gamma}$.
Given   $f,g, h \in \F W(L^p,\ell^q) (\mathbb R^d),$   then  $H_{\gamma}(f, g, h) \in  \F W(L^p,\ell^q) (\R^d)$, and the following estimate  holds
$$\|H_{\gamma}(f, g, h)\|_{\F W(L^p,\ell^q)} \lesssim \|f\|_{\F W(L^p,\ell^q)}  \|g\|_{\F W(L^p,\ell^q)}  \|h\|_{\F W(L^p,\ell^q)}.$$
\end{proposition}
\begin{proof}  By Corollary \ref{m1}, we have
\begin{eqnarray*}
\| H_{\gamma}(f, g,h)\|_{\F W(L^p,\ell^q)} &  \lesssim &    \|  |\cdot|^{-\gamma} \ast (f \bar{g})\|_{\mathcal{F}L^{1}} \|h\|_{\F W(L^p,\ell^q)}.
\end{eqnarray*}  We note that
\begin{eqnarray*}
\left| |\xi|^{-(d-\gamma)} \widehat{f \bar{g}}(\xi) \right|
& \leq &  \tfrac{1}{|\xi|^{d-\gamma}}  \int_{\mathbb R^d} |\widehat{f}(\xi - \eta)| | \widehat{\bar{g}}(\eta)| d\eta
\end{eqnarray*}
and integrating with respect to $\xi,$ we get
\begin{eqnarray*}
  \|  |\cdot|^{-\gamma} \ast (f \bar{g})\|_{\mathcal{F}L^{1}} & \lesssim & \int_{\R^d} \int_{\R^d} \frac{|\widehat{f}(\xi_1)| |\widehat{\bar{g}}(\xi_2)|}{|\xi_1- \xi_2|^{d-\gamma}} d\xi_1 d\xi_2= \left\langle | I^{\gamma}\widehat{f}|,| \widehat{\bar{g}}|\right \rangle_{L^2(\R^d)}
\end{eqnarray*}
where $I^{\gamma}$ denotes the Riesz potential of order $\gamma$:
\[
I^{\gamma}\widehat{f}(x)=C_{\gamma}\int_{\mathbb R^d}\tfrac{\widehat{f}(y)}{|x-y|^{d-\gamma}} dy.\] By H\"older and  Hardy-Littlewood Sobolev  inequalities and Lemma \ref{pl}\eqref{in}, we have
\begin{eqnarray*}
\|  |\cdot|^{-\gamma} \ast (f \bar{g})\|_{\mathcal{F}L^{1}} &= & \|I^{\gamma}\widehat{f}\|_{L^{\frac{2d}{d-\gamma}}} \|\widehat{\bar{g}}\|_{L^{\frac{2d}{d+\gamma}}}\\
  & \lesssim  & \|\widehat{f}\|_{L^{\frac{2d}{d+\gamma}}} \|\widehat{\bar{g}}\|_{L^{\frac{2d}{d+\gamma}}}= \|f\|_{\F W(L^{\frac{2d}{d+\gamma}},\ell^{ \frac{2d}{d+\gamma}})} \|g\|_{\F W(L^{\frac{2d}{d+\gamma}},\ell^{\frac{2d}{d+\gamma}})}\\
  & \lesssim &  \|f\|_{\F W(L^p,\ell^q)} \|g\|_{\F W(L^p,\ell^q)}.
\end{eqnarray*}
This completes the proof.
\end{proof}

\begin{proposition}\label{t4}
Let $0<\gamma<d$ and $f, g, h \in \F W(L^p,\ell^q)\cap L^2(\R^d).$
\begin{enumerate}
\item \label{s1}  Let  $1\leq q<\frac{2d}{d+\gamma}(<2)$ and $q\leq p$. Then
$$\|H_{\gamma}(f, g, h)\|_{\F W(L^p,\ell^q)\cap L^2} \lesssim \|f\|_{\F W(L^p,\ell^q)\cap L^2}  \|g\|_{\F W(L^p,\ell^q)\cap L^2}  \|h\|_{\F W(L^p,\ell^q)\cap L^2}.$$
\item \label{s2}  Let  $(2<)\frac{2d}{d-2\gamma} <q\leq\infty$ and $p\leq q$. %and $0<\gamma<d(\frac{1}{2}-\frac{1}{q})$.
Then
$$\|H_{\gamma}(f, g, h)\|_{\F W(L^p,\ell^q)} \lesssim \|f\|_{\F W(L^p,\ell^q)\cap L^2}  \|g\|_{\F W(L^p,\ell^q)\cap L^2}  \|h\|_{\F W(L^p,\ell^q)\cap L^2}.$$
\end{enumerate}
\end{proposition}
\begin{proof}
By Corollary \ref{m1} and \eqref{k1k2},  we have
\begin{eqnarray}
\|H_{\gamma}(f, g, h)\|_{\F W(L^p,\ell^q)}
&\lesssim &\|  |\cdot|^{-d+\gamma}  (\widehat{f}\ast\widehat{ \bar{g}})\|_{L^1}\|h\|_{\F W(L^p,\ell^q)}\nonumber\\
&\lesssim & \left(\| k_1 (\widehat{f}\ast\widehat{ \bar{g}})\|_{L^1}+\|  k_2  (\widehat{f}\ast\widehat{ \bar{g}})\|_{L^1}\right)\|h\|_{\F W(L^p,\ell^q)}. \label{n3}
\end{eqnarray}
By H\"older and Hausdorff-Young inequalities,  we have
\begin{eqnarray}\label{n4}
 \| k_1 (\widehat{f}\ast\widehat{ \bar{g}})\|_{L^1}\leq \| k_1\|_{L^1}\| \widehat{f}\ast\widehat{ \bar{g}}\|_{L^\infty}\lesssim \| \widehat{f\bar{g}}\|_{L^\infty}
 \leq \| f\bar{g}\|_{L^1}
\leq \| f\|_{L^2}\|g\|_{L^2}
\end{eqnarray}
\eqref{s1} Note that
$$\tfrac{1}{q/[2(q-1)]}+\tfrac{1}{q/[2-q]}=1, \quad \tfrac{1}{q}+\tfrac{1}{q}=1+\tfrac{1}{q/[2-q]}, $$
and in view of  \eqref{k1k2},  we impose the following condition:  $$\tfrac{q}{2(q-1)}>\tfrac{d}{d-\gamma}\Longleftrightarrow q<\tfrac{2d}{d+\gamma}.$$ Thus,  by H\"older and Young inequalities, for $q\leq p,$  we obtain
\begin{eqnarray}\label{n5}
\|  k_2  (\widehat{f}\ast\widehat{ \bar{g}})\|_{L^1}&\leq&\|k_2\|_{{q}/{(2(q-1))}}\|\widehat{f}\ast\widehat{ \bar{g}}\|_{q/(2-q)}\nonumber\\
&\lesssim&\|\widehat{f}\|_{L^q}\|\widehat{\bar{g}}\|_{L^q}=\|f\|_{\F L^q}\|g\|_{\F L^q}\nonumber\\
& \lesssim &  \|f\|_{\F W(L^p,\ell^q)} \|g\|_{\F W(L^p,\ell^q)}.
\end{eqnarray}
Combining  \eqref{n3}, \eqref{n4} and \eqref{n5},  we obtain
\[
\|H_{\gamma}(f, g, h)\|_{\F W(L^p,\ell^q)} \lesssim \left(\| f\|_2\|g\|_2 + \|f\|_{\F W(L^p,\ell^q)} \|g\|_{\F W(L^p,\ell^q)}\right)\|h\|_{\F W(L^p,\ell^q)}
\]
Taking $p=q=2$ in \eqref{n3}   we have
\begin{eqnarray*}
\|H_{\gamma}(f, g, h)\|_{L^2} &\lesssim& \left(\| k_1 (\widehat{f}\ast\widehat{ \bar{g}})\|_{L^1}+\|  k_2  (\widehat{f}\ast\widehat{ \bar{g}})\|_{L^1}\right)\|h\|_2\\
&\lesssim&\left(\| f\|_2\|g\|_2+\|f\|_{\F W(L^p,\ell^q)} \|g\|_{\F W(L^p,\ell^q)}\right)\|h\|_{L^2}.
\end{eqnarray*} using \eqref{n4}, \eqref{n5}. This completes proof of \eqref{s1}.\\\\
\eqref{s2} In view of \eqref{k1k2},
we may rewrite
\begin{eqnarray}\label{n7}
H_{\gamma}(f,g,h)=\left( k_1^\vee \ast (f \bar{g}) \right) h+\left( k_2^\vee \ast (f \bar{g}) \right) h,
\end{eqnarray}
where $k_i^{\vee}$ denotes the inverse Fourier transform of $k_i.$ By Corollary \ref{m1} and \eqref{n4},  we obtain
\begin{eqnarray}\label{n8}
\|\left( k_1^\vee \ast (f \bar{g}) \right) h\|_{\F W(L^p,\ell^q)}&\lesssim&\|k_1( \widehat{f}\ast\widehat{ \bar{g}}\|_{L^1}\| h\|_{\F W(L^p,\ell^q)}\lesssim  \| f\|_{L^2}\|g\|_{L^2}\| h\|_{\F W(L^p,\ell^q)}
\end{eqnarray}
Note that
$$\tfrac{1}{2q/(q+2)}+\tfrac{1}{2}=1+\tfrac{1}{q},$$
and in view of \eqref{k1k2} we impose condition:
$\frac{2q}{q+2}>\frac{d}{d-\gamma}\Leftrightarrow q>\frac{2d}{d-2\gamma}$.
Since $p\leq q$,  by Lemma \ref{pl}\eqref{in}, Corollary \ref{m1},  and Lemma \ref{pl} \eqref{con},    we obtain
\begin{eqnarray}\label{n6}
\|\left( k_2^\vee \ast (f \bar{g}) \right) h\|_{\F W(L^p,\ell^q)}&\le &\|\left( k_2^\vee \ast (f \bar{g}) \right) h\|_{\F W(L^q,\ell^q)}\\
&\leq&\| k_2^\vee \ast (f \bar{g})\|_{\widehat{w}^{2q/(q+2),2q/(q+2)}}\|h\|_{L^2}\nonumber\\
&\leq&\| k_2^\vee  \|_{\F W(L^{2q/(q+2)},\ell^{2q/(q+2)})}\|f\bar{g}\|_{\F L^\infty}\|h\|_2\nonumber\\
&\lesssim&\|\widehat{f\bar{g}}\|_{\infty}\|h\|_2\nonumber\\
&\lesssim&\|{f\bar{g}}\|_{1}\|h\|_{L^2}
\leq\|f\|_{L^2}\|g\|_2\|h\|_{L^2}
\end{eqnarray}
Now  using \eqref{n7}, \eqref{n8} and \eqref{n6},  we have
\begin{eqnarray*}
\|H_{\gamma}(f, g, h)\|_{\F W(L^p,\ell^q)}&\lesssim& \| f\|_{L^2}\|g\|_{L^2}\| h\|_{\F W(L^p,\ell^q)}
+\|f\|_{L^2}\|g\|_{L^2}\|h\|_{L^2}
\end{eqnarray*}which completes the proof.
\end{proof}

\begin{remark} The trilinear estimate in Proposition \ref{t4} \eqref{s2} is not invoked in this article, this may be of  independent interest.
\end{remark}

\begin{proposition}\label{t3}
Let $0<\gamma<\frac{d}{2}$, $1\leq p,q\leq\infty$,  $\frac{d}{d-\gamma}< \rho\leq2$ and $Y= \F W(L^p,\ell^q)$ or $M^{p,q}.$ Then
$$\|H_{\gamma}(f, g, h)\|_{Y} \lesssim(\left\| f\right\|_{L^2}\left\| g\right\|_{L^2}+\left\| f\right\|_{L^{2\rho}}\left\| f\right\|_{L^{2\rho}})  \left\| h\right\|_{Y}.$$
\end{proposition}
\begin{proof} By \eqref{k1k2} and Corollary \ref{m1} and Proposition \ref{des}\eqref{mfl}, we obtain
\begin{align*}
\left\|\h_{\gamma}(f,g,h)\right\|_{Y}&\lesssim\|K\ast(f\bar{g})\|_{\F L^1}\|h\|_{Y}
=\left\|\F K\F(f\overline{g})\right\|_{L^1}\left\| h\right\|_{Y}\\
& \leq \left(\left\| k_1\right\|_{L^1}\left\|\F(f\overline{g})\right\|_{L^\infty}+\left\| k_2\right\|_{L^\rho}\left\|\F(f\overline{g})\right\|_{L^{\rho'}} \right) \left\| h\right\|_{Y}\\
&\lesssim (\left\| f\overline{g}\right\|_{L^1}+\left\| f\overline{g}\right\|_{L^\rho})  \left\| h\right\|_{Y}\\
& \leq  (\left\| f\right\|_{L^2}\left\| g\right\|_{L^2}+\left\| f\right\|_{L^{2\rho}}\left\| f\right\|_{L^{2\rho}})  \left\| h\right\|_{Y}.\qedhere
\end{align*}
\end{proof}

\begin{definition}\label{fpd} Let $s\in [0,1]$. Any pair $(q,r)$ of positive real number is said to be $s-$admissible, if $q, r\geq 2$ and
\[ \tfrac{2s}{q}+ \tfrac{d}{r}= \tfrac{d}{2}.\]
Such set of all $s-$admissible pair is denoted by $\Gamma_s$.
\end{definition}
\begin{theorem}[Stichartz estimates]\label{stri} \
\begin{enumerate}
\item\cite[Theorem 1]{chergui2022blowup},  \cite[Corollary 3.10]{guo2014improved} \label{stri1} Let $d\geq 2,  \frac{d}{2d-1}< s_1\leq s_2\leq 1$ and $u_0, u, F$ are radial in space and satisfying \eqref{mfh}. Then
\[\|u\|_{L^q_tL^r_x} \lesssim \|u_0\|_{2}+\|F\|_{L^{\tilde{q}'}_t L^{\tilde{r}'}_x} \]
if $(q,r)$ and $(\tilde{q}, \tilde{r})$ belong to $\Gamma_{s_1} \cup \Gamma_{s_2}$ and either $(\tilde{q}', \tilde{r}')\neq (2, \infty)$ or $(q,r)\neq (2, \infty)$.
\item\cite{keel1998endpoint} \label{stri2} For $s_1=s_2=1,$ the above  estimate holds  for any $d\geq 1$ without any radiality assumption.
\end{enumerate}
 \end{theorem}
 \subsection{Harmonic Oscillator}
\begin{theorem}(uniform boundedness of linear propogator, \cite[Theorem 5]{drt}, cf. \cite{cordero2008metaplectic})\label{mso}
The Schr\"odinger  propagator  associated to the harmonic oscillator $e ^{itH}$ is bounded on $M^{p,p}(\mathbb R^d)$ for each $t\in \mathbb R$, and all  $1\leq p < \infty.$ Moreover, we have $\|e^{itH}f\|_{M^{p,p}}=\|f\|_{M^{p,p}}.$
\end{theorem}
%\begin{definition} A pair $(q,r)$ is admissible if  $2\leq r< \frac{2d}{d-2}$ with $2\leq r \leq \infty$ if $d=1$, and  $2\leq r < \infty$ if $d=2$, whenever
%$$\frac{2}{q} =  d \left( \frac{1}{2} - \frac{1}{r} \right).$$
%\end{definition}

\begin{proposition}(Strichartz estimates, \cite[Proposition 2.2]{carles2011nonlinear}) \label{seh}    Let $\phi \in L^2(\mathbb R^d)$ and
$$DF(t,x) :=  (e^{itH}\phi)(x) +  \int_0^t e^{i(t-\tau )H}F(\tau,\cdot)(x) d\tau.$$  Then for any time slab $I$ and  $1-$admissible pairs $(q_j,r_j)$, $j=1,2$ with $q_j>2$,
there exists  a constant $C=C(|I|, r_1)$ such that for all intervals $I \ni 0, $
$$ \|D(F)\|_{L_t^{q_1}L_x^{r_1}}  \leq  C \|\phi \|_{L^2}+   C  \|F\|_{L_t^{q_2'}L_x^{r_2'}}, \ \forall \phi\in L^2,\forall F \in L^{q_2'} (I, L^{r_2'}).$$
\end{proposition}
\section{Proof of main results}
\subsection{Local wellposedness}
As mentioned earlier, the proof uses a fixed point argument.  We have divided the proof in two cases.  In the Case I,  we use the trilinear estimates in Propositions \ref{t2}, and \ref{t4},  whereas in Case II, Proposition  \ref{t3} along with Strichartz estimates (Proposition \ref{des}) is used.
\begin{proof}[{\bf Proof of Theorem \ref{MT}}]
{\bf Case I:}
$
X=\begin{cases} \F W(L^p,\ell^q) & \text{with}  \    1 \leq q \leq \frac{2d}{d+\gamma}\le p \leq \infty\\
\F W(L^p,\ell^q)\cap L^2 & \text{with} \    1\leq q<\frac{2d}{d+\gamma}, q\leq p\leq\infty
\end{cases}.
$
By Duhamel's formula, we write \eqref{mfh}
as \begin{equation}\label{m}
u(t)= U(t)u_{0}-i\int_{0}^{t}U(t-\tau) \, [(K\ast |u|^{2}(\tau))u(\tau)] \, d\tau:=\mathcal{J}_{u_0}(u).
\end{equation}
By Lemma \ref{pl}\eqref{n1} and Propositions \ref{t2} and \ref{t4}\eqref{s1},  for $u_0\in B_m$ we have
\begin{eqnarray}\label{tac}
\|\mathcal{J}_{u_0}u\|_{C([0, T], X)} \leq  C \left(\|u_{0}\|_{X} + c T \|u\|_{X}^{3}\right)\leq  C \left(m + c T \|u\|_{X}^{3}\right),
\end{eqnarray}
for some universal constant $c.$

For $M>0$, set  $U_{T, M}= \{u\in C([0, T],  X):\|u\|_{C([0, T],  X)}\leq M \}.$
We claim that  $\mathcal{J}_{u_0}: U_{T, M} \to U_{T, M}$ is a contraction mapping for a suitable choice of  $M$ and small $T>0$.  Indeed for $u_0\in B_m$, if we let, $M= 2Cm$ and $u\in U_{T, M},$ from \eqref{tac} we obtain
\begin{eqnarray}
\|\mathcal{J}_{u_0}u\|_{C([0, T],  X)} \leq  \tfrac{M}{2} + cC_{T}T M^{3}.
\end{eqnarray}
We choose a  $T$  such that  $c CTM^{2} \leq 1/2,$ that is, $T=T(m, d) \sim m^{-2} $ and as a consequence  we have
\begin{eqnarray}
\|\mathcal{J}_{u_0}u\|_{C([0, T],  X)} \leq \tfrac{M}{2} + \tfrac{M}{2}=M,
\end{eqnarray}
therefore, $\mathcal{J}_{u_0}u \in U_{T, M}.$
By the arguments as before and using trilinearity of $H_\gamma$, for $u_0,v_0\in B_m$ and $u,v\in U_{T,M}$ we obtain
\begin{eqnarray}\label{n2}
\|\mathcal{J}_{u_0}u- \mathcal{J}_{v_0}v\|_{C([0, T],  X)} \leq \|u_0-v_0\|_X+\tfrac{1}{2} \|u-v\|_{C([0, T],  X)}.
\end{eqnarray}
Therefore, (putting $v_0=u_0$) using the  Banach's contraction mapping principle, we conclude that $\mathcal{J}_{u_0}$ has a fixed point say $u$ in $U_{T, M}$ which is a solution of \eqref{mfh}. Also if $v$ is the fixed point of $\mathcal{J}_{v_0}$ in $U_{T, M}$, from \eqref{n2}, it follows that
\[
\|u- v\|_{C([0, T],  X)} \leq \|u_0-v_0\|_X+\tfrac{1}{2} \|u-v\|_{C([0, T],  X)}\Longrightarrow \|u- v\|_{C([0, T],  X)} \leq 2\|u_0-v_0\|_X
\]
which shows the solution map is Lipschitz continuous.\\

\textbf{Case II:} $X=\begin{cases}\F W(L^p,\ell^q)\cap L^2 &  \text{with}\  d\geq 1,  s_1= s_2 = 1, \gamma<\min(2,\frac{d}{2})\\
\F W(L^p,\ell^q)\cap L_{rad}^2 &  \text{with} \     d\geq 2,  \frac{d}{2d-1}< s_1\leq s_2 \leq 1, \gamma<\min(2s_2,\frac{d}{2})
\end{cases}$\\
 or \\
\hspace*{1.7cm}$
X=\begin{cases}
M^{p,q}\cap L^2 &  \text{with}\  d\geq 1,  s_1= s_2 = 1, \gamma<\min(2,\frac{d}{2})\\
M^{p,q}_{rad}\cap L^2 &  \text{with} \     d\geq 2,  \frac{d}{2d-1}< s_1\leq s_2 \leq 1, \gamma<\min(2s_2,\frac{d}{2})
\end{cases}$.\\

Let  ${q_1}=\frac{8s_2}{\gamma}, r=\frac{4d}{2d-\gamma}.$ Then $(q_1,r)$ is $s_2-$admissible ($q_1\geq2$ is ensured as $\gamma<2s_2< 4s_2$,  see Definition \ref{fpd}).
For  $T, b>0,$  introduce the space
$$U_b^T=\big\{v\in L_T^\infty\big(X\big):\| v\|_{L_T^\infty(X)}\leq b,\| v\|_{L_T^{q_1}(L^r)}\leq b, \| v\|_{L_T^{2q_2}(L^{2\rho})}\leq b\big\},$$
where $q_2,\rho$ to be chosen later.   We set %$\mathcal{U}_b^T=(U_b^T)^N$ and
  the  distance on it by
$$d(u,v)=\max\big\{\| u-v\|_{L_T^\infty(X)},\| u-v\|_{L_T^{q_1}(L^r)},\| u-v\|_{L_T^{2q_2}(L^{2\rho})}\big\},$$
where $u,v\in U_b^T$.  Next, we show that the mapping $\mathcal{J}_{u_0},$ defined by  \eqref{m}, takes ${U}_b^T$ into itself for suitable choice of $b$ and small  $T>0$.  Let $u,v,w\in {U}_b^T$ and
denote
\begin{align}\label{f1}
J(t)=\int_0^tU(t-s)H_\gamma(u(s),v(s),w(s))ds.
\end{align}
Let $\frac{d}{d-\gamma}<\rho\leq2,$
and choose  $q_2>1$  as
 $\frac{2s_2}{2q_2} =  d \big( \frac{1}{2} - \frac{1}{2\rho} \big)$
 so that $(2q_2,2\rho)$ is an $s_2-$fractional admissible pair ($2q_2>2$ imposes the condition $\gamma<2s_2$).
 Then using Proposition \ref{t3} with $Y=\F W(L^p,\ell^q)$ or $M^{p,q}$ (radiality condition is incorporated in definition of $Y$ if $s_1\neq1$),
\begin{align*}
\left\| J(t)\right\|_{Y}&\lesssim\int_0^t\big(\left\|u(s)\right\|_{L^2}\left\|v(s)\right\|_{L^2}+\left\|u(s)\right\|_{L^{2\rho}}\left\|v(s)\right\|_{L^{2\rho}}\big)\left\|w(s)\right\|_{Y}ds\\
&\lesssim t\left\|u\right\|_{L_t^\infty(L^2)}\left\|v\right\|_{L_t^\infty(L^2)}\left\|w\right\|_{L_t^\infty(Y)}+\left\|u\right\|_{L_t^{2q_2}\left(L^{2\rho}\right)}\left\|v\right\|_{L_t^{2q_2}\left(L^{2\rho}\right)}\left\|w\right\|_{L_t^{{q_2}'}\left(Y\right)}\nonumber
\end{align*}
using H\"older inequality.
Therefore by this and Lemma \ref{pl}\eqref{n1}, Proposition \ref{des}\eqref{des2}, for $u_0\in B_m$
\begin{align}\label{12a}
\left\|\mathcal{J}_{u_0}(u)(t)\right\|_{Y}\lesssim\left\|u_0\right\|_{Y}+b^3(T+T^{\frac{1}{q_2'}})\le m+2b^3T^{1-\frac{1}{q_2}}.
\end{align}
For $(\underline{q},\underline{r})\in\{({q_1},r),(2q_2,2\rho),(\infty,2)\},$ by Proposition \ref{stri}  we have
\begin{align*}
\|\mathcal{J}_{u_0}(u)\|_{L^{\underline{q}}([0,T),L^{\underline{r}})}&\lesssim \|u_0\|_{L^2}+ \| (K*|u|^2)u\|_{L^{{q_1}'}([0,T),L^{r'})}.
\end{align*}
Now we have
 $\frac{1}{{q_1}'}=\frac{4s_2-\gamma}{4s_2}+\frac{1}{{q_1}},\ \frac{1}{r'}=\frac{\gamma}{2d}+\frac{1}{r}\ \text{and}\ \frac{4s_2-\gamma}{4s_2}=\frac{2}{{q_1}}+\frac{2s_2-\gamma}{2s_2}.$
By H\"older and Hardy-Littlewood-Sobolev inequalities, %as in the calculation \eqref{b1}
we have
\begin{eqnarray}\label{n9}
\left\| (K*(u\overline{v})w\right\|_{L^{{q_1}'}([0,T),L^{r'})}
&\leq&\bigl\|\bigl\||\cdot|^{-\gamma}*|u\overline{v}|\bigr\|_{L^{\frac{2d}{\gamma}}}\bigr\|_{L^{\frac{4s_2}{4s_2-\gamma}}([0,T))}\left\|w\right\|_{L^{q_1}([0,T),L^r)}\nonumber\\
&\lesssim&\bigl\|\left\||u\overline{v}|\right\|_{L^{\frac{2d}{2d-\gamma}}}\bigr\|_{L^{\frac{4s_2}{4s_2-\gamma}}([0,T))}\left\|w\right\|_{L^{q_1}([0,T),L^r)}\nonumber\\
&\leq& T^{1-\frac{\gamma}{2s_2}}\left\|u\right\|_{L^{q_1}([0,T),L^r)}\left\|v\right\|_{L^{q_1}([0,T),L^r)}\left\|w\right\|_{L^{q_1}([0,T),L^r)}.
\end{eqnarray}
Combining the above two inequalities, we obtain
\begin{align*}%\label{7b}
\|\mathcal{J}_{u_0}(u)\|_{L^{\underline{q}}(L^{\underline{r}})}&\lesssim\|u_0\|_{L^2}+ T^{1-\frac{\gamma}{2s_2}}\left\|u\right\|_{L_T^{q_1}(L^r)}^3\lesssim m+ T^{1-\frac{1}{q_2}}b^3.
\end{align*}
 as $1-\frac{1}{q_2}<1-\frac{\gamma}{2s_2}\Leftrightarrow\frac{d}{d-\gamma}<\rho$.
 Choose $b=2cm$ and $T\sim m^{-\frac{2q_2}{q_2-1}}>0$ small enough so that    \eqref{12a} and the above inequality %\eqref{7b}
  imply  $\mathcal{J}_{u_0}(u) \in {U}_b^T$.
On the other hand for $u_0,v_0\in B_m$, $u,v\in {U}_b^T,$ using trilinearity of $H_{\gamma}$, proceeding as above, we have
\begin{eqnarray}\label{12b}
&&\|\mathcal{J}_{u_0}(u)(t)-\mathcal{J}_{v_0}(v)(t)\|_{Y}\nonumber\\
&\leq&\|u_0-v_0\|_{Y}+\int_0^t\left\|H_{\gamma}(u(s),u(s),u(s))-H_{\gamma}(v(s),v(s),v(s))\right\|_{Y}ds\nonumber\\
&\lesssim& \|u_0-v_0\|_{Y}+T^{1-\frac{1}{q_2}}b^2d(u,v)
\end{eqnarray}
and
\begin{align}\label{7d}
\| \mathcal{J}_{u_0}(u)-\mathcal{J}_{v_0}(v)\|_{L^{\underline{q}}(L^{\underline{r}})}\lesssim \|u_0-v_0\|_{L^2}+T^{1-\frac{1}{q_2}}Nb^2d(u,v).
\end{align}
Choose $T>0$ further small so that   \eqref{12b}  and \eqref{7d} (after putting $u_0=v_0$) imply that $\mathcal{J}_{u_0}$ is a contraction.
Continuous dependence of solution on data follows from \eqref{12b}, \eqref{7d} as in the first case.
%In view of Propositions \ref{des} and \ref{t3} and exploiting the above arguments,  the last case of Theorem \ref{MT} follows.  Thus,  we shall omit the details.
\end{proof}
\begin{proof}[\textbf{Proof of Theorem \ref{gHh}} \eqref{gHh1}]
 Taking Proposition \ref{seh} and Theorem \ref{mso} into account,  the proof of Theorem \ref{gHh} \eqref{gHh1} follows as the Case II of the above proof.

 By Duhamel's formula, we write \eqref{hH}
as $u(t)=\mathcal{I}_{u_0}(u)$ where $\mathcal{I}_{u_0}(u)$ is defined by
\begin{equation*}
\mathcal{I}_{u_0}(u):= e^{itH}u_{0}-i\int_{0}^{t}e^{i(t-\tau)H} \, [(K\ast |u|^{2}(\tau))u(\tau)] \, d\tau.
\end{equation*}
 Let  ${q_1}=\frac{8}{\gamma}, r=\frac{4d}{2d-\gamma}$ so that $(q_1,r)$ becomes an admissible pair and set
 $$U_b^T=\big\{v\in L_T^\infty(M^{p,p}\cap L^2):\| v\|_{L_T^\infty( M^{p,p}\cap L^2)}\leq b,\| v\|_{L_T^{q_1}(L^r)}\leq b, \| v\|_{L_T^{2q_2}(L^{2\rho})}\leq b\big\},$$
 and distance on it
 $$d(u,v)=\max\big\{\| u-v\|_{L_T^\infty(M^{p,p}\cap L^2)},\| u-v\|_{L_T^{q_1}(L^r)},\| u-v\|_{L_T^{2q_2}(L^{2\rho})}\big\},$$
 with $q_2,\rho$ to be chosen such a way that
 $\frac{d}{d-\gamma}<\rho\leq2,$
and  $q_2>1$  as
 $\frac{2}{2q_2} =  d \big( \frac{1}{2} - \frac{1}{2\rho} \big)$
 so that $(2q_2,2\rho)$ is an $1-$fractional admissible pair ($2q_2>2$ imposes the condition $\gamma<2$).
 For $u,v,w\in {U}_b^T$,
denote
\begin{align*}
I(t)=\int_0^te^{i(t-s)H}H_\gamma(u(s),v(s),w(s))ds.
\end{align*}
By Proposition \ref{t3} with $Y=M^{p,p}$
\begin{align*}
\left\| I(t)\right\|_{M^{p,p}}&\lesssim\int_0^t\big(\left\|u(s)\right\|_{L^2}\left\|v(s)\right\|_{L^2}+\left\|u(s)\right\|_{L^{2\rho}}\left\|v(s)\right\|_{L^{2\rho}}\big)\left\|w(s)\right\|_{M^{p,p}}ds\\
&\lesssim t\left\|u\right\|_{L_t^\infty(L^2)}\left\|v\right\|_{L_t^\infty(L^2)}\left\|w\right\|_{L_t^\infty(M^{p,p})}+\left\|u\right\|_{L_t^{2q_2}\left(L^{2\rho}\right)}\left\|v\right\|_{L_t^{2q_2}\left(L^{2\rho}\right)}\left\|w\right\|_{L_t^{{q_2}'}\left(M^{p,p}\right)}\nonumber
\end{align*}
using H\"older inequality.
Then%refore by this and
Theorem \ref{mso}, %Proposition \ref{seh},
for $u_0\in B_m$
\begin{align}\label{h1}
\left\|\mathcal{J}_{u_0}(u)(t)\right\|_{M^{p,p}}\lesssim\left\|u_0\right\|_{M^{p,p}}+b^3(T+T^{\frac{1}{q_2'}})\le m+2b^3T^{1-\frac{1}{q_2}}.
\end{align}
For $(\underline{q},\underline{r})\in\{({q_1},r),(2q_2,2\rho),(\infty,2)\},$ by Proposition \ref{seh}  we have
\begin{align*}
\|\mathcal{I}_{u_0}(u)\|_{L^{\underline{q}}([0,T),L^{\underline{r}})}&\lesssim \|u_0\|_{L^2}+ \| (K*|u|^2)u\|_{L^{{q_1}'}([0,T),L^{r'})}.
\end{align*}
Note that
 $\frac{1}{{q_1}'}=\frac{4-\gamma}{4}+\frac{1}{{q_1}},\ \frac{1}{r'}=\frac{\gamma}{2d}+\frac{1}{r}\ \text{and}\ \frac{4-\gamma}{4}=\frac{2}{{q_1}}+\frac{2-\gamma}{2}.$
By H\"older and Hardy-Littlewood-Sobolev inequalities, %as in the calculation \eqref{b1}
we have
\begin{eqnarray*}
\left\| (K*(u\overline{v})w\right\|_{L^{{q_1}'}([0,T),L^{r'})}
&\leq&\bigl\|\left\||\cdot|^{-\gamma}*|u\overline{v}|\right\|_{L^{\frac{2d}{\gamma}}}\bigr\|_{L^{\frac{4}{4-\gamma}}([0,T))}\left\|w\right\|_{L^{q_1}([0,T),L^r)}\nonumber\\
&\lesssim&\bigl\|\left\||u\overline{v}|\right\|_{L^{\frac{2d}{2d-\gamma}}}\bigr\|_{L^{\frac{4}{4-\gamma}}([0,T))}\left\|w\right\|_{L^{q_1}([0,T),L^r)}\nonumber\\
&\leq& T^{1-\frac{\gamma}{2}}\left\|u\right\|_{L^{q_1}([0,T),L^r)}\left\|v\right\|_{L^{q_1}([0,T),L^r)}\left\|w\right\|_{L^{q_1}([0,T),L^r)}.
\end{eqnarray*}
Combining the above two inequalities, we obtain
\begin{align*}%\label{7b}
\|\mathcal{J}_{u_0}(u)\|_{L^{\underline{q}}(L^{\underline{r}})}&\lesssim\|u_0\|_{L^2}+ T^{1-\frac{\gamma}{2}}\left\|u\right\|_{L_T^{q_1}(L^r)}^3\lesssim m+ T^{1-\frac{1}{q_2}}b^3.
\end{align*}
 as $1-\frac{1}{q_2}<1-\frac{\gamma}{2}\Leftrightarrow\frac{d}{d-\gamma}<\rho$.
 Choose $b=2cm$ and $T\sim m^{-\frac{2q_2}{q_2-1}}>0$ small enough so that    \eqref{h1} and the above inequality %\eqref{7b}
  imply  $\mathcal{J}_{u_0}(u) \in {U}_b^T$.
On the other hand for $u_0,v_0\in B_m$, $u,v\in {U}_b^T,$ using trilinearity of $H_{\gamma}$, proceeding as above, we have
\begin{eqnarray*}
&&\|\mathcal{I}_{u_0}(u)(t)-\mathcal{I}_{v_0}(v)(t)\|_{M^{p,p}}\nonumber\\
&\leq&\|u_0-v_0\|_{M^{p,p}}+\int_0^t\left\|H_{\gamma}(u(s),u(s),u(s))-H_{\gamma}(v(s),v(s),v(s))\right\|_{M^{p,p}}ds\nonumber\\
&\lesssim& \|u_0-v_0\|_{M^{p,p}}+T^{1-\frac{1}{q_2}}b^2d(u,v)
\end{eqnarray*}
and
\begin{align*}
\| \mathcal{I}_{u_0}(u)-\mathcal{I}_{v_0}(v)\|_{L^{\underline{q}}(L^{\underline{r}})}\lesssim \|u_0-v_0\|_{L^2}+T^{1-\frac{1}{q_2}}Nb^2d(u,v).
\end{align*}
Choose $T>0$ further small so that   the above two inequalities %\eqref{12b}  and \eqref{7d} (after putting $u_0=v_0$)
 imply that $\mathcal{I}_{u_0}$ is a contraction.
Continuous dependence of solution on data also follows from these inequalities.% \eqref{12b}, \eqref{7d} as in the first case.
 %We shall omit the details.
\end{proof}
\subsection{Global wellposedness in $L^2$}
Using the estimate we prove the existence of a global solution for $u_0\in L^2$. Moreover these solutions will have additional regularity %. The latter %additional regularity
 which will be instrumental to achieve the global solution from the local one established in Theorem \ref{MT}.
\begin{proposition}[global well-posedness in $L^2$]\label{miD}
Let $d\geq 2,$  $\frac{d}{2d-1} < s_1\leq s_2 \leq 1,$  and  $K$ be given by \eqref{hk} with $\lambda \in \mathbb R$. Fix $j\in\{1,2\}$ and $0<\gamma < \text{min} \{2s_j, d\}.$ If $u_{0}\in L_{rad}^{2}(\mathbb R^{d}),$ then %for each $j\in\{1,2\}$,
 \eqref{mfh} has a unique global solution
$$u\in C(\mathbb R, L_{rad}^{2})\cap L^{8s_2/\gamma}_{loc}(\mathbb R, L^{4d/(2d-\gamma)}).$$
In addition, its $L^{2}-$norm is conserved,
$$\|u(t)\|_{L^{2}}=\|u_{0}\|_{L^{2}}, \   \forall t \in \mathbb R,$$
and for all %$\alpha-$ fractional admissible pairs
$(q,r)\in\Gamma_{s_1} \cup \Gamma_{s_2}$, $u \in L_{loc}^q(\mathbb R, L^r(\mathbb R^d)).$ In the case $s_1=s_2=1$, the radiallity condition on $u_0$ and dimension restriction can be removed.
\end{proposition}
\begin{proof}
% By Duhamel's formula, we write \eqref{mfh}
%as
%$$u(t)=S(t)u_0- i \int_0^t S(t-\tau) (K \ast |u|^2)u(\tau) d\tau:= \Phi(u)(t).$$
 Let us fix a $j\in\{1,2\}$. We introduce the space
\begin{eqnarray*}
Y_j(T)  & =  &\big\{ \phi \in C\left([0,T], L^2 \right): \|\phi \|_{L^{\infty}([0, T], L^2)}  \leq b,\  \|\phi\|_{L^{q_1} ([0,T], L^{r})}  \leq b\big\}
\end{eqnarray*}
 where $q_1= \frac{8s_j}{\gamma}, \  r= \frac{4d}{2d- \gamma}$
and the distance
$$d_j(\phi_1, \phi_2)= \max\left\{ \|\phi_1-\phi_2 \|_{L^{\infty}([0, T], L^2)}, \|\phi_1 - \phi_2 \|_{L^{q_1 }\left( [0, T], L^r\right)} \right\}.$$ Then $(Y_j(T), d_j)$ is a complete metric space. Now we show that $\Phi$ takes $Y_j(T)$ to $Y_j(T)$ for some $T>0.$
%We put
%$$ \ q=q_j= \frac{8s_j}{\gamma}, \  r= \frac{4d}{2d- \gamma}$$
%so that $(q,r)$ is $s_i-$fractional admissible.
Then for $(\underline{q}, \underline{r}) \in \{ (q_1,r), (\infty, 2) \}$,  proceeding as \eqref{n9}, for $\|u_0\|_2\leq m$  we obtain
\begin{equation}\label{n10}
\|\mathcal{J}_{u_0}(u)\|_{L^{\underline{q}}(L^{\underline{r}})}   \lesssim   \|u_0\|_{L^2}+T^{1- \frac{\gamma}{2s_j}} \|u\|^3_{L^{q_1}(L^r)}\le m+T^{1- \frac{\gamma}{2s_j}}b^3
\end{equation}
 where $\mathcal{J}_{u_0}(u)$ is as defined in \eqref{m}.
This shows that $\mathcal{J}_{u_0}$ maps $Y_j(T)$ to $Y_j(T)$ with $b=2c\|u_0\|_{L^2}$ and $T\sim m^{-\frac{4s_j}{2s_j-\gamma}}>0$ small enough.  Next, using trilinearity
$$(K\ast |v|^{2})v- (K\ast |w|^{2})w= (K\ast |v|^{2})(v-w) + (K \ast (|v|^{2}- |w|^{2}))w, $$
 and \eqref{n9} we show $\mathcal{J}_{u_0}$
is a  contraction.
 Then there exists a unique $u \in Y_j(T)$ solving \eqref{mfh}. The global  existence of the solution \eqref{mfh} follows from the conservation of the $L^2-$norm of $u.$
 The last property of the proposition then follows from \eqref{n10} by choosing $(\underline{q},\underline{r})$ %the Strichartz estimates applied with
 an arbitrary $s_1$ or $s_2-$ fractional admissible pair on the left hand side and the same pairs as above on the right hand side.
\end{proof}
Taking Proposition \ref{seh} and Theorem \ref{mso} into account,  we notice that the  the following analogue of Proposition \ref{miD} hold true for \eqref{hH}.
\begin{proposition}[see Proposition 4 in \cite{bhimani2018nonlinear}]\label{h2}
Let  $0<\gamma < \text{min} \{2, d\}.$ If $u_{0}\in L^{2},$ then \eqref{hH} has a unique global solution
$$u\in C(\R,  L^{2})\cap L^{8/\gamma}_{loc}(\R, L^{4d/(2d-\gamma)} ).$$
In addition, its $L^{2}-$norm is conserved,
$$\|u(t)\|_{L^{2}}=\|u_{0}\|_{L^{2}}, \   \forall t \in \mathbb R,$$
and for all  admissible pairs  $(p,q), u \in L_{loc}^{p}(\mathbb R, L^{q}).$
\end{proposition}
\subsection{Global wellposedness in Fourier amalgam spaces and in modulation spaces}
Now we show that with the help of Proposition \ref{miD}, the local solution achieved in Theorem \ref{MT} can be extended on full real line $\R$.
\begin{proof}[\textbf{Proof of Theorem \ref{mt2}}]
We shall see that the solution constructed before in Theorem \ref{MT} is global in time if $u_0\in L^2$ (and radial if $s_1\neq1$).  First note that by redefining
\[
\mathcal{J}_{u_0}(u):= U(t-t_0)u_{0}-i\int_{t_0}^{t}U(t- \tau)\,  \, [(K\ast |u|^{2}(\tau))u(\tau)] \, \, d\tau
\] and proceeding as above in the case $t_0=0$, we can find a solution around time $t_0$ for duration $T=T(\|u(t_0)\|_X,d)$. This shows that the solution can be continued throughout time if  $\|u(t)\|_{X}$ stays bounded in finite time interval. (This proves the blow-up alternative.)

 In fact,
  using Proposition \ref{miD}, we will show that
 $\|u(t)\|_{X}$ cannot become unbounded in finite time.

 Assume $0<T_{\max}$ be so that \eqref{mfh} has a solution in $[0,T_{max})$ and let $0<T<T_{\max}$. %In view of \eqref{dc} and to use the Hausdorff-Young inequality we
Let $1< \frac{d}{d-\gamma} <r \leq 2,$
and using Proposition \ref{t3}, we obtain
\begin{eqnarray}
\|u(t)\|_{X} & \lesssim &  \|u_{0}\|_{X} + \int_{0}^{t} \|H_\gamma(u(\tau),u(\tau),u(\tau))\|_{X}d\tau \nonumber \\
&\lesssim &  \|u_{0}\|_{X} + \int_{0}^{t}(\|u(s)\|_2^2+\|u(\tau)\|_{L^{2r}}^{2})\|u(\tau)\|_{X}d\tau\nonumber\\
%& \lesssim &  C_{T} \left ( \|u_{0}\|_{X} + \int_{0}^{t} \|K\ast |u(\tau)|^{2}\|_{\mathcal{F}L^{1}} \|u(\tau)\|_{X} d\tau \right) \nonumber\\
%& \lesssim & C_{T}  \|u_{0}\|_{X} + C_{T}\int_{0}^{t} \left( \|k_{1}\|_{L^{1}} \|u(\tau)\|_{L^{2}}^{2}+ \|k_{2}\|_{L^{r}} \|\widehat{|u(\tau)|^{2}}\|_{L^{r'}}
%\right) \nonumber \\
%&&\|u(\tau)\|_{X} d\tau \nonumber\\
%& \lesssim & C_{T}\|u_{0}\|_{X} + C_{T} \int_{0}^{t} \left(  \|k_{1}\|_{L^{1}} \|u_{0}\|_{L^{2}}^{2}+ \|k_{2}\|_{L^{r}} \||u(\tau) |^{2}\|_{L^{r}}\right) \nonumber  \\
%&& \|u(\tau)\|_{X} d\tau\nonumber\\
& \lesssim & \|u_{0}\|_{X}+ \int_{0}^{t}\|u(\tau)\|_{X} d\tau  + \int_{0}^{t} \|u(\tau)\|_{L^{2r}}^{2} \|u(\tau)\|_{X } d\tau,\nonumber
 \end{eqnarray}
 where we have used   %H\"older's inequality, and
 the conservation of the $L^{2}-$norm of $u$ from Proposition \ref{miD}.
We note that the requirement on $r$ can be fulfilled if and only if $0<\gamma <d/2.$

%To apply Proposition \ref{mi}, we
Let $\beta>1$ be so that  $(2\beta, 2 r)$ is  $s_2-$fractional admissible, that is, $\frac{s_2}{\beta}= d \left(\frac{1}{2}- \frac{1}{2r} \right)$ such that $\frac{1}{\beta}= \frac{d}{2s_2} \left( 1 - \frac{1}{r} \right)<1.$ This is possible provided that $\frac{r-1}{r} < \frac{2s_2}{d}:$ this condition is compatible with the requirement $r> \frac{d}{d-\gamma}$ if and only if $\gamma < 2s_2.$
 Using the H\"older's inequality for the last integral above, for $0\leq t\leq T$ we obtain
\begin{eqnarray*}
\|u(t)\|_{X} &  \leq & c\|u_0\|_{X} + c\int_{0}^{t} \|u(\tau)\|_{X} d\tau  + c\|u\|_{L^{2\beta}([0, t], L^{2r})}^{2}\|u\|_{L^{\beta'}[(0, t],X)}\\
&  \leq & c\|u_0\|_{X} + c\int_{0}^{t} \|u(\tau)\|_{X} d\tau  + c\|u\|_{L^{\beta'}([0, t],X)},
\end{eqnarray*}
as $\|u\|_{L^{2\beta}([0, T_{\max}), L^{2r})}<\infty$ by Proposition \ref{miD}. % where $\beta'$ is the H\"older conjugate exponent of $\beta.$
Set,
$$h(t):=\sup_{0 \leq \tau \leq t} \|u(\tau)\|_{X}.$$
%For a given $T>0,$
Then $h$ satisfies an estimate of the form,
$$h(t)\leq c\|u_{0}\|_{X}+ c\int_{0}^{t} h(\tau) d\tau + c \left( \int_{0}^{t}h(\tau)^{\beta'} d\tau \right)^{\frac{1}{\beta'}},$$
provided that $0 \leq t \leq T,$ and where we have used the fact that $\beta'$ is finite.
Using the H\"older's inequality we infer that,
$$h(t)\leq  c \|u_{0}\|_{X} + C(T) \left(\int_{0}^{t} h(\tau)^{\beta'}d \tau \right)^{\frac{1}{\beta'}},\qquad C(T)=c(T^{\frac{1}{\beta}}+1).$$
Raising the above estimate to the power $\beta'$, we find that
$$h(t)^{\beta'} \leq  C_{2}(T,\|u_0\|_X) \left( 1+\int_{0}^{t} h(\tau)^{\beta'} d\tau\right).$$
In view of  Gr\"onwall's inequality, we conclude that  $h\in L^{\infty}([0, T]).$ Since $0<T<T_{\max}$ is arbitrary, $h\in L^{\infty}_{loc}([0,T_{\max})),$ and  the proof of Theorem \ref{mt2}  follows.
\end{proof}
\begin{proof}[\textbf{Proof of Theorem \ref{gHh} \eqref{gH2}}]
%Taking Remarks  \ref{DS1} and \ref{DS2} into account and exploiting  the method of the  proof of Theorem \ref{MT},  the proof follows.

%\textcolor{red}{Can you write a sketch of proof here?}
As in the above proof, it is enough to show $\|u(t)\|_{M^{p,p}\cap L^2}$  cannot become unbounded in finite time. Assume $0<T_{\max}$ be so that \eqref{hH} has a solution in $[0,T_{max})$, let $0<T<T_{\max}$ and $r$ be as in the above proof. Then using Proposition \ref{t3} with $Y=M^{p,p}$
\begin{eqnarray}
\|u(t)\|_{M^{p,p}\cap L^2}  \lesssim
\|u_{0}\|_{M^{p,p}\cap L^2}+ \int_{0}^{t}\|u(\tau)\|_{M^{p,p}\cap L^2} d\tau  + \int_{0}^{t} \|u(\tau)\|_{L^{2r}}^{2} \|u(\tau)\|_{M^{p,p}\cap L^2 } d\tau,\nonumber
 \end{eqnarray} using conservation of the $L^{2}-$norm of $u$ from Proposition \ref{h2}.

Let $\beta>1$ be so that  $(2\beta, 2 r)$ is  $1-$fractional admissible, that is, $\frac{1}{\beta}= d \left(\frac{1}{2}- \frac{1}{2r} \right)$ such that $\frac{1}{\beta}= \frac{d}{2s_2} \left( 1 - \frac{1}{r} \right)<1.$ This is possible provided that %$\frac{r-1}{r} < \frac{2s_2}{d}:$ this condition is compatible with the requirement $r> \frac{d}{d-\gamma}$ if and only if
$\gamma < 2.$
Then we have $0\leq t\leq T$ we obtain 
\begin{eqnarray*}
\|u(t)\|_{M^{p,p}\cap L^2} %&  \leq & c\|u_0\|_{M^{p,p}\cap L^2} + c\int_{0}^{t} \|u(\tau)\|_{M^{p,p}\cap L^2} d\tau  + c\|u\|_{L^{2\beta}([0, t], L^{2r})}^{2}\|u\|_{L^{\beta'}[(0, t],M^{p,p}\cap L^2)}\\
&  \leq & c\|u_0\|_{M^{p,p}\cap L^2} + c\int_{0}^{t} \|u(\tau)\|_{M^{p,p}\cap L^2} d\tau  + c\|u\|_{L^{\beta'}([0, t],M^{p,p}\cap L^2)},
\end{eqnarray*}
as $\|u\|_{L^{2\beta}([0, T_{\max}), L^{2r})}<\infty$ by Proposition \ref{h2}. % where $\beta'$ is the H\"older conjugate exponent of $\beta.$
Set,
$$h(t):=\sup_{0 \leq \tau \leq t} \|u(\tau)\|_{M^{p,p}\cap L^2}.$$
%For a given $T>0,$
and proceed as in the above proof by replacing $X$ with $M^{p,p}\cap L^2$.
\end{proof}

\noindent
{\textbf{Acknowledgement}:} S Haque is thankful to DST--INSPIRE  (DST/INSPIRE/04/2022/001457) \& USIEF--Fulbright-Nehru fellowship for financial support. S Haque is also thankful to Harish-Chandra Research Institute \&  University of California, Los Angeles  for their excellent research facilities. 
\bibliographystyle{siam}
\bibliography{mixed}

\begin{thebibliography}{10}

\bibitem{bejenaru2006sharp}
{\sc I.~Bejenaru and T.~Tao}, {\em
  \href{https://doi.org/10.1016/j.jfa.2005.08.004}{Sharp well-posedness and
  ill-posedness results for a quadratic non-linear {S}chr\"{o}dinger
  equation}}, J. Funct. Anal., 233 (2006), pp.~228--259.

\bibitem{benyi2009local}
{\sc A.~B\'{e}nyi and K.~A. Okoudjou}, {\em Local well-posedness of nonlinear
  dispersive equations on modulation spaces}, Bull. Lond. Math. Soc., 41
  (2009), pp.~549--558.

\bibitem{kassob}
{\sc {\'A}.~B{\'e}nyi and K.~A. Okoudjou}, {\em {M}odulation {S}paces: {W}ith
  {A}pplications to {P}seudodifferential {O}perators and {N}onlinear
  {S}chr\"odinger Equations,}, 2020.

\bibitem{bhimani2018nonlinear}
{\sc D.~G. Bhimani}, {\em The nonlinear {S}chr{\"o}dinger equations with
  harmonic potential in modulation spaces}, Discrete \& Continuous Dynamical
  Systems-A, 39 (2019), pp.~5923-- 5944.

\bibitem{bhimani2023nodea}
{\sc D.~G. Bhimani}, {\em The blow-up solutions for fractional heat equations
  on torus and {E}uclidean space}, NoDEA Nonlinear Differential Equations
  Appl., 30 (2023), p.~Paper No. 19.

\bibitem{drt}
{\sc D.~G. Bhimani, R.~Balhara, and S.~Thangavelu}, {\em Hermite multipliers on
  modulation spaces}, in In: Delgado J., Ruzhansky M. (eds) Analysis and
  Partial Differential Equations: Perspectives from Developing Countries.
  Springer Proceedings in Mathematics {\&} Statistics, vol 275, Springer, Cham,
  2019.

\bibitem{bhimani2020hartree}
{\sc D.~G. Bhimani, M.~Grillakis, and K.~A. Okoudjou}, {\em The
  {H}artree--{F}ock equations in modulation spaces}, Communications in Partial
  Differential Equations,  (2020), pp.~1--30.

\bibitem{D3}
{\sc D.~G. Bhimani, M.~Grillakis, and K.~A. Okoudjou}, {\em The
  {H}artree-{F}ock equations in modulation spaces}, Comm. Partial Differential
  Equations, 45 (2020), pp.~1088--1117.

\bibitem{bhimani2023sharp}
{\sc D.~G. Bhimani, H.~Hajaiej, S.~Haque, and T.~Luo}, {\em A sharp
  gagliardo-nirenberg inequality and its application to fractional problems
  with inhomogeneous nonlinearity}, Evolution Equations and Control Theory, 12
  (2023), pp.~362--390.

\bibitem{bhimani2023hartree}
{\sc D.~G. Bhimani and S.~Haque}, {\em
  \href{https://doi.org/10.1007/s00023-022-01234-5}{The {H}artree and
  {H}artree-{F}ock equations in {L}ebesgue {$L^p$} and {F}ourier-{L}ebesgue
  {$\hat{L}^p$} spaces}}, Ann. Henri Poincar\'{e}, 24 (2023), pp.~1005--1049.

\bibitem{bhimaniAdM2021}
{\sc D.~G. Bhimani, R.~Manna, F.~Nicola, S.~Thangavelu, and S.~I. Trapasso},
  {\em Phase space analysis of the {H}ermite semigroup and applications to
  nonlinear global well-posedness}, Adv. Math., 392 (2021), pp.~Paper No.
  107995, 18.

\bibitem{bhimani2022heat}
\leavevmode\vrule height 2pt depth -1.6pt width 23pt, {\em On heat equations
  associated with fractional harmonic oscillators}, Fract. Calc. Appl. Anal.,
  26 (2023), pp.~2470--2492.

\bibitem{bhimani2016functions}
{\sc D.~G. Bhimani and P.~K. Ratnakumar}, {\em Functions operating on
  modulation spaces and nonlinear dispersive equations}, J. Funct. Anal., 270
  (2016), pp.~621--648.

\bibitem{biagi2022mixed}
{\sc S.~Biagi, S.~Dipierro, E.~Valdinoci, and E.~Vecchi}, {\em Mixed local and
  nonlocal elliptic operators: regularity and maximum principles}, Comm.
  Partial Differential Equations, 47 (2022), pp.~585--629.

\bibitem{bourgain1999JAMS}
{\sc J.~Bourgain}, {\em Global wellposedness of defocusing critical nonlinear
  {S}chr\"{o}dinger equation in the radial case}, J. Amer. Math. Soc., 12
  (1999), pp.~145--171.

\bibitem{carles2011nonlinear}
{\sc R.~Carles}, {\em Nonlinear {S}chr{\"o}dinger equation with time dependent
  potential}, Communications in Mathematical Sciences, 9 (2011), pp.~937--964.

\bibitem{carles2014cauchy}
{\sc R.~Carles and L.~Mouzaoui}, {\em
  \href{https://doi.org/10.1090/S0002-9939-2014-12072-7 }{On the {C}auchy
  problem for the {H}artree type equation in the {W}iener algebra}},
  Proceedings of the American Mathematical Society, 142 (2014), pp.~2469--2482.

\bibitem{leonid2017JDE}
{\sc L.~Chaichenets, D.~Hundertmark, P.~Kunstmann, and N.~Pattakos}, {\em On
  the existence of global solutions of the one-dimensional cubic {NLS} for
  initial data in the modulation space {$M_{p,q}(\Bbb R)$}}, J. Differential
  Equations, 263 (2017), pp.~4429--4441.

\bibitem{chergui2022blowup}
{\sc L.~Chergui}, {\em On blowup solutions for the mixed fractional
  {S}chr\"{o}dinger equation of {C}hoquard type}, Nonlinear Anal., 224 (2022),
  pp.~Paper No. 113105, 21.

\bibitem{cho2017AA}
{\sc Y.~Cho, M.~M. Fall, H.~Hajaiej, P.~A. Markowich, and S.~Trabelsi}, {\em
  Orbital stability of standing waves of a class of fractional
  {S}chr\"{o}dinger equations with {H}artree-type nonlinearity}, Anal. Appl.
  (Singap.), 15 (2017), pp.~699--729.

\bibitem{hichem2013FE}
{\sc Y.~Cho, H.~Hajaiej, G.~Hwang, and T.~Ozawa}, {\em On the {C}auchy problem
  of fractional {S}chr\"{o}dinger equation with {H}artree type nonlinearity},
  Funkcial. Ekvac., 56 (2013), pp.~193--224.

\bibitem{cho2006SIAM}
{\sc Y.~Cho and T.~Ozawa}, {\em On the semirelativistic {H}artree-type
  equation}, SIAM J. Math. Anal., 38 (2006), pp.~1060--1074.

\bibitem{cordero2008metaplectic}
{\sc E.~Cordero and F.~Nicola}, {\em Metaplectic representation on {W}iener
  amalgam spaces and applications to the {S}chr{\"o}dinger equation}, Journal
  of Functional Analysis, 254 (2008), pp.~506--534.

\bibitem{deng2013estimate}
{\sc Q.~Deng, Y.~Ding, and L.~Sun}, {\em Estimate for generalized unimodular
  multipliers on modulation spaces}, Nonlinear {A}nalysis: {T}heory, {M}ethods
  {\&} {A}pplications, 85 (2013), pp.~78--92.

\bibitem{FeiPAMS}
{\sc H.~Feichtinger and A.~Gumber}, {\em Completeness of shifted dilates in
  invariant {B}anach spaces of tempered distributions}, Proc. Amer. Math. Soc.,
  149 (2021), pp.~5195--5210.

\bibitem{Fei}
{\sc H.~G. Feichtinger}, {\em Banach convolution algebras of {W}iener type}, in
  Functions, series, operators, {V}ol. {I}, {II} ({B}udapest, 1980), vol.~35 of
  Colloq. Math. Soc. J\'{a}nos Bolyai, North-Holland, Amsterdam, 1983,
  pp.~509--524.

\bibitem{feichtinger1983modulation}
{\sc H.~G. Feichtinger}, {\em Modulation spaces on locally compact {A}belian
  groups}.
\newblock Technical Report, University of Vienna, 1983, and in “Wavelets and
  Their Applications” (eds. M. Krishna, R. Radha and S. Thangavelu), 99-140,
  Allied Publishers, New Delhi, 2003., 1983.
\newblock Available on
  %\href{https://www.researchgate.net/profile/Hans_Feichtinger/publication/200524248_Modulation_spaces_on_locally_compact_Abelian_group/links/542002980cf2218008d43064.pdf}
  {researchgate.net}.

\bibitem{HGF2006}
\leavevmode\vrule height 2pt depth -1.6pt width 23pt, {\em Modulation spaces:
  looking back and ahead}, Sampl. Theory Signal Image Process., 5 (2006),
  pp.~109--140.

\bibitem{HF2015}
\leavevmode\vrule height 2pt depth -1.6pt width 23pt, {\em Choosing function
  spaces in harmonic analysis}, in Excursions in harmonic analysis. {V}ol. 4,
  Appl. Numer. Harmon. Anal., Birkh\"{a}user/Springer, Cham, 2015, pp.~65--101.

\bibitem{JH}
{\sc J.~Forlano and T.~Oh}, {\em Normal form approach to the one-dimensional
  cubic nonlinear {S}chr{\"o}dinger equation in {F}ourier-amalgam spaces,
  \textit{preprint}}.

\bibitem{forlano2020deterministic}
{\sc J.~A. Forlano}, {\em On the deterministic and probabilistic {C}auchy
  problem of nonlinear dispersive partial differential equations}, PhD thesis,
  Heriot-Watt University, 2020.

\bibitem{grochenig2013foundations}
{\sc K.~Gr\"{o}chenig}, {\em Foundations of time-frequency analysis}, Applied
  and Numerical Harmonic Analysis, Birkh\"{a}user Boston, Inc., Boston, MA,
  2001.

\bibitem{GuoJDE}
{\sc Q.~Guo and S.~Zhu}, {\em Sharp threshold of blow-up and scattering for the
  fractional {H}artree equation}, J. Differential Equations, 264 (2018),
  pp.~2802--2832.

\bibitem{guo2014improved}
{\sc Z.~Guo and Y.~Wang}, {\em Improved {S}trichartz estimates for a class of
  dispersive equations in the radial case and their applications to nonlinear
  {S}chr{\"o}dinger and wave equations}, Journal d'Analyse Math{\'e}matique,
  124 (2014), pp.~1--38.

\bibitem{hajaiej2014multilinear}
{\sc H.~Hajaiej, P.~A. Markowich, and S.~Trabelsi}, {\em Multiconfiguration
  {H}artree-{F}ock theory for pseudorelativistic systems: the time-dependent
  case}, Math. Models Methods Appl. Sci., 24 (2014), pp.~599--626.

\bibitem{hajaiej2011necessary}
{\sc H.~Hajaiej, L.~Molinet, T.~Ozawa, and B.~Wang}, {\em Necessary and
  sufficient conditions for the fractional {G}agliardo-{N}irenberg inequalities
  and applications to {N}avier-{S}tokes and generalized boson equations}, in
  Harmonic analysis and nonlinear partial differential equations, RIMS
  K\^{o}ky\^{u}roku Bessatsu, B26, Res. Inst. Math. Sci. (RIMS), Kyoto, 2011,
  pp.~159--175.

\bibitem{CH2002}
{\sc C.~Heil}, {\em An introduction to weighted {W}iener amalgams}, Wavelets
  and their Applications (Chennai, January 2002), Krishna, M., Radha, R., and
  Thangavelu, S., Eds., AlliEd., Publishers, New Delhi, 183–216.

\bibitem{CH2007}
\leavevmode\vrule height 2pt depth -1.6pt width 23pt, {\em History and
  evolution of the density theorem for {G}abor frames}, J. Fourier Anal. Appl.,
  13 (2007), pp.~113--166.

\bibitem{HerrNA2014}
{\sc S.~Herr and E.~Lenzmann}, {\em The {B}oson star equation with initial data
  of low regularity}, Nonlinear Anal., 97 (2014), pp.~125--137.

\bibitem{HyakunaMulti}
{\sc R.~Hyakuna}, {\em
  \href{https://doi.org/10.1007/s00028-018-0432-8}{Multilinear estimates with
  applications to nonlinear {S}chr\"{o}dinger and {H}artree equations in
  {$\widehat{L^p}$}-spaces}}, J. Evol. Equ., 18 (2018), pp.~1069--1084.

\bibitem{hyakuna2012PAMS}
{\sc R.~Hyakuna and M.~Tsutsumi}, {\em On existence of global solutions of
  {S}chr\"{o}dinger equations with subcritical nonlinearity for
  {$\widehat{L}^p$}-initial data}, Proc. Amer. Math. Soc., 140 (2012),
  pp.~3905--3920.

\bibitem{keel1998endpoint}
{\sc M.~Keel and T.~Tao}, {\em Endpoint {S}trichartz estimates}, American
  Journal of Mathematics, 120 (1998), pp.~955--980.

\bibitem{klausJFAA2023}
{\sc F.~Klaus}, {\em Wellposedness of {NLS} in modulation spaces}, J. Fourier
  Anal. Appl., 29 (2023), pp.~Paper No. 9, 37.

\bibitem{magin2008modeling}
{\sc R.~Magin and M.~Ovadia}, {\em Modeling the cardiac tissue electrode
  interface using fractional calculus}, Journal of Vibration and Control, 14
  (2008), pp.~1431--1442.

\bibitem{manna2022JDE}
{\sc R.~Manna}, {\em On the existence of global solutions of the {H}artree
  equation for initial data in the modulation space {$M^{p,q}(\Bbb R)$}}, J.
  Differential Equations, 317 (2022), pp.~70--88.

\bibitem{changxing2008cauchy}
{\sc C.~Miao, G.~Xu, and L.~Zhao}, {\em The {C}auchy problem of the {H}artree
  equation}, Journal of Partial Differential Equations, 21 (2008), pp.~22--44.

\bibitem{ruzhansky2012modulation}
{\sc M.~Ruzhansky, M.~Sugimoto, and B.~X. Wang}, {\em Modulation spaces and
  nonlinear evolution equations}, in Evolution equations of hyperbolic and
  {S}chr\"{o}dinger type, vol.~301 of Progr. Math., Birkh\"{a}user/Springer
  Basel AG, Basel, 2012, pp.~267--283.

\bibitem{vargas2001JMPA}
{\sc A.~Vargas and L.~Vega}, {\em Global wellposedness for 1{D} non-linear
  {S}chr\"{o}dinger equation for data with an infinite {$L^2$} norm}, J. Math.
  Pures Appl. (9), 80 (2001), pp.~1029--1044.

\bibitem{wang2007global}
{\sc B.~Wang and H.~Hudzik}, {\em The global {C}auchy problem for the {NLS} and
  {NLKG} with small rough data}, J. Differential Equations, 232 (2007),
  pp.~36--73.

\bibitem{wang2011harmonic}
{\sc B.~Wang, Z.~Huo, C.~Hao, and Z.~Guo}, {\em
  \href{https://doi.org/10.1142/9789814360746}{Harmonic analysis method for
  nonlinear evolution equations. {I}}}, World Scientific Publishing Co. Pte.
  Ltd., Hackensack, NJ, 2011.

\bibitem{wangJFA2006}
{\sc B.~Wang, Z.~Lifeng, and G.~Boling}, {\em Isometric decomposition
  operators, function spaces {$E^\lambda_{p,q}$} and applications to nonlinear
  evolution equations}, J. Funct. Anal., 233 (2006), pp.~1--39.

\bibitem{NW26}
{\sc N.~Wiener}, {\em On the representation of functions by trigonometrical
  integrals}, Math. Z., 24 (1926), pp.~575--616.

\bibitem{NW}
\leavevmode\vrule height 2pt depth -1.6pt width 23pt, {\em Tauberian theorems},
  Ann. of Math. (2), 33 (1932), pp.~1--100.

\end{thebibliography}
\end{document}